\documentclass[journal]{IEEEtran}

\ifCLASSINFOpdf
\else
\fi

\usepackage[final]{graphicx}
\usepackage[cmex10]{amsmath}
\usepackage{amssymb}
\usepackage{stdmath}
\usepackage[lined,boxed,commentsnumbered,vlined,ruled]{algorithm2e}
\usepackage{slashbox,booktabs,amsmath}
\usepackage{subfigure}

\begin{document}
\title{Solving Large-Scale Robust Stability Problems by Exploiting the Parallel Structure of Polya's Theorem}

\author{Reza~Kamyar, Matthew~M.~Peet and Yulia~Peet

\thanks{Manuscript received November 21, 2011; revised September 24, 2012; accepted February 26, 2013. This work was supported entirely by funding from the National Science Foundation under Award CMMI-1100376. Recommended
by XXXX.}
\thanks{R. Kamyar is with the Department of Mechanical Engineering, Cybernetic Systems and Controls Laboratory, Arizona State University, Tempe, AZ 85281 USA (e-mail: {rkamyar@asu.edu \tt\small }).}%
\thanks{M. M. Peet and Y. Peet are with the School for Engineering of Matter, Transport, and Energy, Engineering Research Center, Arizona State University, Tempe, AZ 85281 USA (e-mail: {mpeet@asu.edu; ypeet@asu.edu). \tt\small}}%
\thanks{Color versions of one or more of the figures in this paper are available online at http://ieeexplore.ieee.org.}
\thanks{Digital Object Identifier 10.1109/TAC.2013.2253253}}


\maketitle

\begin{abstract}
In this paper, we propose a distributed computing approach to solving large-scale robust stability problems on the simplex. Our approach is to formulate the robust stability problem as an optimization problem with polynomial variables and polynomial inequality constraints. We  use Polya's theorem to convert the polynomial optimization problem to a set of highly structured Linear Matrix Inequalities (LMIs). We then use a slight modification of a common interior-point primal-dual algorithm to solve the structured LMI constraints. This yields a set of extremely large yet structured computations. 
 We then map the structure of the computations to a decentralized computing environment consisting of independent processing nodes with a structured adjacency matrix. The result is an algorithm which can solve the robust stability problem with the same per-core complexity as the deterministic stability problem with a conservatism which is only a function of the number of processors available. Numerical tests on cluster computers and supercomputers demonstrate the ability of the algorithm to efficiently utilize hundreds and potentially thousands of processors and analyze systems with 100+ dimensional state-space. The proposed algorithms can be extended to perform stability analysis of nonlinear systems and robust controller synthesis.
 
\end{abstract}



\section{Introduction}

\IEEEPARstart{T}{his} paper addresses the problem of stability of large-scale systems with several unknown parameters.
Control system theory when applied in practical situations often involves the use of large state-space models, typically due to inherent complexity of the system, the interconnection of subsystems, or the reduction of an infinite-dimensional or PDE model to a finite-dimensional approximation. One approach to dealing with such large scale models has been to use model reduction techniques such as balanced truncation~\cite{gugercin2004survey}. However, the use of model reduction techniques are not necessarily robust and can result in arbitrarily large errors. In addition to large state-space, practical problems often contain uncertainty in the model due to modeling errors, linearization, or fluctuation in the operating conditions. The problem of stability and control of systems with uncertainty has been widely studied. See, e.g. the texts~\cite{ackermann_2001,bhattacharyya_1995,green_1994,zhou_1998}. However, a limitation of existing computational methods for analysis and control of systems with uncertainty is high complexity. This is a consequence of fact that the problem of robust stability of systems with parametric uncertainty is known to be NP-hard~\cite{blondel2000survey,nemirovskii1993several}. The result is that for systems with parametric uncertainty and with hundreds of states, existing algorithms will fail with the primary point of failure usually being lack of unallocated memory. 

In this paper, we seek to distribute the computation laterally over an array of processors within the context of existing computational resources. Specifically, we seek to utilize cluster-computing, supercomputing and Graphics Processing Unit (GPU)-computing architectures. When designing algorithms to run in a parallel computing environment, one must both synchronize computational tasks among the processors while minimizing communication overhead among the processors. This can be difficult, as each architecture has a specific communication graph. we account for communication by explicitly modeling the required communication graph between processors. This communication graph is then mapped to the processor architecture using the Message-Passing Interface (MPI)~\cite{walker1996mpi}. While there are many algorithms for robust stability analysis and control of linear systems, ours is the first which explicitly accounts for the processing architecture in the emerging multi-core computing environment.

Our approach to robust stability is based on the well-established use of parameter-dependent Quadratic-In-The-State (QITS) Lyapunov functions. The use of parameter-dependent Lyapunov QITS functions eliminates the conservativity associated with e.g. quadratic stability~\cite{quadratic1}, at the cost of requiring some restriction on the rate of parameter variation. Specifically, our QITS Lyapunov variables are polynomials in the vector of uncertain parameters. This is a generalization of the use of QITS Lyapunov functions with affine parameter dependence as in~\cite{affine_orig} and expanded in, e.g.~\cite{affine_dependent0,affine_dependent1,affine_dependent2,affine_dependent3}. The use of polynomial QITS Lyapunov variables can be motivated by~\cite{Bliman_existence}, wherein it is shown that any feasible parameter-dependent LMI with parameters inside a compact set has a polynomial solution or~\cite{peet2009exponentially} wherein it is shown that local stability of a nonlinear vector field implies the existence of a polynomial Lyapunov function.

There are several results which use polynomial QITS Lyapunov functions to prove robust stability. In most cases, the stability problem is reduced to the general problem of optimization of polynomial variables subject to LMI constraints - an NP-hard problem~\cite{Np_hard}. To avoid NP-hardness, the polynomial optimization problem is usually solved in an asymptotic manner by posing a sequence of sufficient conditions of increasing accuracy and decreasing conservatism. For example, building on the result in~\cite{Bliman_existence},~\cite{Bliman_convex} provides a sequence of increasingly precise LMI conditions for robust stability analysis of linear systems with affine dependency on uncertain parameters on the complex unit ball. Necessary and sufficient stability conditions for linear systems with one uncertain parameter are derived in~\cite{bound1}, providing an explicit bound on the degree of the polynomial-type Lyapunov function. The result is extended to multi-parameter-dependent linear systems in~\cite{bound2}. Another important approach to optimization of polynomials is the Sum of Squares (SOS) methodology which replaces the polynomial positivity constraint with the constraint that the polynomial admits a representation as a sum of squares of polynomials~\cite{sos1,sos2,sos3,Homog_chesi}. A version of this theorem for polynomials with matrix coefficients can be found in~\cite{sos3}. While we have worked extensively with the SOS methodology, we have not, as of yet, been able to adapt algorithms for solving the resulting LMI conditions to a parallel-computing environment. Finally, there have been several results in recent years on the use of Polya's theorem to solve polynomial optimization problems~\cite{polya_peres} on the simplex. An extension of the Polya's theorem for uncertain parameters on the multisimplex or hypercube can be found in~\cite{multisimplex}. The approach presented in this paper is an extension of the use of Polya's theorem for solving polynomial optimization problems in a parallel computing environment.

The goal of this project is to create algorithms which explicitly map computation, communication and storage to existing parallel processing architectures. This goal is motivated by the failure of existing general-purpose Semi-Definite Programming (SDP) solvers to efficiently utilize platforms for large-scale computation. Specifically, it is well-established that linear programming and semi-definite programming both belong to the complexity class P-Complete, also known as the class of inherently sequential problems.
Although there have been several attempts to map certain SDP solvers to a parallel computing environment~\cite{csdp,sdpara}, certain critical steps cannot be distributed. The result is that as the number of processors increases, certain bottleneck computations dominate, leading a saturation in computational speed of these solvers (Amdahl's law~\cite{amdahl}). We avoid these bottleneck computations and communications by exploiting the particular structure of the LMI conditions associated with Polya's theorem. Note that, in principle, a perfectly designed general-purpose SDP algorithm could identify the structure of the SDP, as we have, and map the communication, computation and memory constraints to the parallel architecture. Indeed, there has been a great deal of research on creating programming languages which attempt to do just this~\cite{kale1994charm,deitz2005high}. However, at present such languages are mostly theoretical and have certainly not been incorporated into existing SDP solvers.

In addition to parallel SDP solvers, there have been some efforts to exploit structure in certain polynomial optimization algorithms to reducing the size and complexity of the resulting LMI's. For example, in~\cite{parrilo_sym} symmetry was used to reduce the size of the SDP variables. Specific sparsity structure was used in~\cite{parrilo_sructure,kim2005generalized,waki} to reduce the complexity of the linear algebra calculations. Generalized approaches to the use of sparsity in SDP algorithms can be found in~\cite{kim2005generalized}. Groebner basis techniques~\cite{ideals,grobner} have been used by~\cite{parrilo_sructure} to simplify the formulation of the SDPs associated with the SOS decomposition problems.


The paper is organized around two independent problems: setting up the sequence of structured SDPs associated with Polya's theorem and solving them. Note that the problem of decentralizing the set-up algorithm is significant in that for large-scale systems, the instantiation of the problem may be beyond the memory and computational capacity of a single processing node. For the set-up problem, the algorithm that we propose has no centralized memory or computational requirements whatsoever. Furthermore, if a sufficient number of processors are available, the number of messages does not change with the size of the state-space or the number of Polya's iterations. In addition, the ideal communication architecture for the set-up algorithm does not correspond to the communication structure of GPU computing or supercomputing.
In the second problem, we propose a variant of a standard SDP primal-dual algorithm and map the computational, memory and communication requirements to a parallel computing environment. Unlike the set-up algorithm, the primal-dual algorithm does have a small centralized component corresponding to the update of the set of dual variables. However, we have structured the algorithm so that the size of this dual computation is solely a function of the degree of the polynomial QITS Lyapunov function and does not depend on the number of Polya's iterations, meaning that the sequence of algorithms has fixed centralized computational and communication complexity. In addition, there is no communication between processors, which means that the algorithm is well suited to most parallel computing architectures. A graph representation of the communication architecture of both the set-up and SDP algorithms has also been provided in the relevant sections.

Combining the set-up and SDP components and testing the result of both in cluster computing environments, we demonstrate the capability of robust analysis and control of systems with 100+ states and several uncertain parameters. Specifically, we ran a series of numerical experiments using a local Linux cluster and the Blue Gene supercomputer (with 200 processor allocation). First, we applied the algorithm to a current problem in robust stability analysis of magnetic confinement fusion using a discretized PDE model. Next, we examine the accuracy of the algorithm as Polya's iterations progress and compare this accuracy with the SOS approach. We show that unlike the general-purpose parallel SDP solver SDPARA~\cite{sdpara}, the speed-up - the increase in processing speed per additional processor - of our algorithm shows no evidence of saturation. Finally, we calculate the envelope of the algorithm on the Linux cluster in terms of the maximum state-space dimension, number of processors and Polya's iterations.\vspace*{0.1in}

\noindent \textit{NOTATION}

We represent $l-$variate monomials as $\alpha^{{\gamma}} = \prod_{i=1}^l\alpha_i^{\gamma_i}$, where $\alpha \in \mathbb{R}^l $ is the vector of variables and $\gamma \in \mathbb{N}^l$ is the vector of exponents and $\sum_{i=1}^l \gamma_i = d$ is the degree of the monomial. We define $W_d:= \left\lbrace \gamma \in \mathbb{N}^l: \sum_{i=1}^l \gamma_i = d \right\rbrace$ as the totally ordered set of the exponents of $l-$variate monomials of degree $d$, where the ordering is lexicographic. In lexicographical ordering $\gamma \in W_d$ precedes $\eta \in W_d$, if the left most non-zero entry of $\gamma-\eta$ is positive. The lexicographical index of every $ \gamma \in W_d$ can be calculated using the map $\langle{{\cdot}}\rangle:\mathbb{N}^l \rightarrow \mathbb{N}$ defined as~\cite{peet_acc} 
\begin{equation}
\langle{\gamma}\rangle = \sum_{j=1}^{l-1} \sum_{i=1}^{\gamma_i} f\bbl(l-j,d+1-\sum_{k=1}^{j-1}\gamma_k -i\bbr) + 1, \label{eq:lex} 
\end{equation}
where as in~\cite{scheinerman} 
\begin{equation}
f(l,d) :=
\begin{cases}
\hspace*{0.2in}	0 &\text{for} \;\; l = 0 \\
\dbinom{l+d-1}{l-1}= \dfrac{(d+l-1)!}{d!(l-1)!} & \text{for} \;\; l > 0, \label{eq:f}
\end{cases}
\end{equation}
is the cardinality of $W_d$, i.e., the number of $l-$variate monomials of degree $d$. For convenience, we also define the index of a monomial $\alpha^\gamma$ to be $\langle \gamma \rangle$. We represent $l-$variate homogeneous polynomials of degree $d_p$ as 
\begin{equation}
P(\alpha)=\sum_{\gamma \in W_{d_p}} P_{\langle{\gamma}\rangle} \alpha^{\gamma}, 
\end{equation}
where $P_{\langle{\gamma}\rangle} \in \mathbb{R}^{n \times n}$ is the matrix coefficient of the monomial $\alpha^{\gamma}$. We denote the element corresponding to the $i^{th}$ row and $j^{th}$ column of matrix $A$ as $[A]_{i,j}$. The subspace of symmetric matrices in $\mathbb{R}^{n \times n} $ is denoted by $ \mathbb{S}^n$.  We define a basis for $\mathbb{S}^n$ as 
\begin{equation*}
[E_k]_{i,j}:=\begin{cases}
1  &  \text{if } i=j=k \\
0  &  \text{otherwise}
\end{cases} , \quad \text{for} \;\, k \leq n \quad
\text{and}  
\end{equation*}
\begin{equation}
[E_k]_{i,j}:=[F_k]_{i,j}+[F_k]^T_{i,j}, \quad \text{for} \;\, k > n,  
\end{equation}
where 
\begin{equation}
[F_k]_{i,j}:=\begin{cases}
1  &  \text{if }  i=j-1=k-n  \\
0  &  \text{otherwise}.
\end{cases} 
\end{equation}
Note that this choice of basis is arbitrary - any other basis could be used. However, any change in basis would require modifications to the formulae defined in this paper.
The canonical basis for $ \mathbb{R}^n $ is denoted by $e_i$ for $i=1, \cdots, n$, where $ 
e_i=[0 \: ... \: 0 \overbrace{1}^{i^{th}} 0 \: ... \: 0].
$
The vector with all entries equal to one is denoted by $\vec{1}$. The trace of $A \in \mathbb{R}^{n \times n}$ is denoted by $tr(A)=\sum_{i=1}^{n} [A]_{i,i}$. The block-diagonal matrix with diagonal blocks $X_1, \cdots, X_m \in \mathbb{R}^{n \times n}$ is denoted $\text{diag}(X_1, \cdots, X_m) \in \mathbb{R}^{mn \times mn} $ or occasionally as $\text{diag}(X_i|_{i=1}^{m}) \in \mathbb{R}^{mn \times mn}$. The identity and zero matrices are denoted by $I_n \in \mathbb{R}^{n \times n}$ and $0_n \in \mathbb{R}^{n \times n}$.


\section{PRELIMINARIES}
\label{sec:PRELIM} 

Consider the linear system 
\begin{equation} 
 \dot{x}(t)= A(\alpha)x(t),  \label{system}
\end{equation}
where $ A(\alpha) \in \mathbb{R}^{n \times n} $ and $\alpha \in Q \subset \R^{l}$ is a vector of uncertain parameters.
In this paper, we consider the case where $A(\alpha)$ is a homogeneous polynomial and
$Q=\Delta_l \subset \R^{l}$ where $\Delta_l$ is the unit simplex, i.e., 
\begin{equation}
\Delta_l=\left\lbrace \alpha\in \mathbb{R}^l , \sum_{i=1}^{l} \alpha_i=1, \alpha_i\geqslant 0 \right\rbrace.  \label{eq:simplex} 
\end{equation}
If $A(\alpha)$ is not homogeneous, we can obtain an equivalent homogeneous representation in the following manner.
Suppose $A(\alpha)$ is a non-homogeneous polynomial with $\alpha \in \Delta_l$, is of degree $d_a$ and has $N_a$ monomials with non-zero coefficients. Define $D = \left( d_{a_1}, \cdots, d_{a_{N_a}} \right)$, where $d_{a_i}$ is the degree of the $i^{th}$ monomial of $A(\alpha)$ according to lexicographical ordering.
Now define the polynomial $B(\alpha)$ as per the following.
\begin{enumerate}
\item Let $B=A$.
\item For $i=1,\cdots, N_a$, multiply the $i^{th}$ monomial of $B(\alpha)$, according to lexicographical ordering, by $\left( \sum_{j=1}^l \alpha_j \right)^{d_a-d_{a_i}}$.
\end{enumerate}
Then, since  $\sum_{j=1}^l \alpha_j =1$,  $B(\alpha)=A(\alpha)$ for all $\alpha \in \Delta_l$ and hence all properties of $\dot x(t) = A(\alpha)x(t)$ are retained by the homogeneous system $\dot x(t) = B(\alpha)x(t)$. \vspace*{0.05in}

\textit{1) Example:} \textit{Construction of the homogeneous system $\dot x(t) = B(\alpha) x(t)$.}

Consider the non-homogeneous polynomial $
A(\alpha) = C \alpha_1^2  + D \alpha_2 + E \alpha_3 + F
$ of degree $d_a=2$, where $[\alpha_1,\alpha_2,\alpha_3] \in \Delta_3$. Using the above procedure, the homogeneous polynomial $B(\alpha)$ can be constructed as 
\begin{align}
&B(\alpha) = C \alpha_1^2  + D \alpha_2 (\alpha_1+\alpha_2+\alpha_3) + E  \alpha_3 (\alpha_1+\alpha_2+\alpha_3) \nonumber \\ 
&+ F (\alpha_1+\alpha_2+\alpha_3)^2  = (\underbrace{C+F}_{B_1})\alpha_1^2 + \underbrace{(D+2F)}_{B_2} \alpha_1 \alpha_2 \nonumber \\ 
&+ \underbrace{(E+2F)}_{B_3} \alpha_1 \alpha_3 + \underbrace{(D+F)}_{B_4} \alpha_2^2 + \underbrace{(D+E+2F)}_{B_5} \alpha_2 \alpha_3 \nonumber \\
&+ \underbrace{(E+F)}_{B_6} \alpha_3^2 = \sum_{\gamma \in W_2} B_{\langle{\gamma}\rangle} \alpha^{\gamma}.
\end{align}

The following is a stability condition~\cite{polya_peres} for System~\eqref{system}.
\begin{thm}\label{thm:thm1}
System~\eqref{system} is stable if and only if there exists a polynomial matrix $P(\alpha)$ such that $P(\alpha) \succ 0$ and
\begin{equation}
A^T(\alpha)P(\alpha)+P(\alpha)A(\alpha) \prec 0   \label{eq:Lyap_LMI}  
\end{equation}
for all $\alpha \in \Delta_l$.
\end{thm}
A similar condition also holds for discrete-time linear systems. The conditions associated with Theorem~\ref{thm:thm1} are infinite-dimensional LMIs, meaning they must hold at infinite number of points. Such problems are known to be NP-hard~\cite{Np_hard}. In this paper we derive a sequence of polynomial-time algorithms such that their outputs converge to the solution of the infinite-dimensional LMI. Key to this result is Polya's Theorem~\cite{polya_book}. A variation of this theorem for matrices is given as follows.

\begin{thm}
\label{Polya's Theorem}
\emph{(Polya's Theorem)}
The homogeneous polynomial $ F(\alpha) \succ 0 $ for all $ \alpha \in\Delta_l $ if and only if for all sufficiently large $d$, 
\begin{equation}
\left(\sum_{i=1}^l \alpha_i \right)^d F(\alpha)  \label{eq:polya's exponent}
\end{equation}
has all positive definite coefficients.
\end{thm}

Upper bounds for Polya's exponent $d$ can be found as in~\cite{d_bound}. However, these bounds are based on the properties of $F$ and are difficult to determine a priori. In this paper, we show that applying Polya's Theorem to the robust stability problem, i.e., the inequalities in Theorem~\ref{thm:thm1} yields a semi-definite programming condition with an efficiently distributable structure. This is discussed in the following section.


\section{PROBLEM SET-UP}
\label{sec:SETUP} 
In this section, we show how Polya's theorem can be used to determine the robust stability of an uncertain system using linear matrix inequalities with a distributable structure.

\subsection {Polya's Algorithm}
We consider the stability of the system described by Equation~\eqref{system}. We are interested in finding a $P(\alpha)$ which satisfies the conditions of Theorem~\ref{thm:thm1}. According to Polya's theorem, the constraints of Theorem~\ref{thm:thm1} are satisfied if for some sufficiently large $d_1$ and $d_2$, the polynomials  
\begin{equation}
\left(\sum_{i=1}^l \alpha_i \right)^{d_1} P(\alpha ) \label{eq:LMI_1} \qquad \text{and} 
\end{equation}
\begin{equation}
-\left(\sum_{i=1}^l \alpha_i \right)^{d_2} \left(A^T(\alpha)P(\alpha)+P(\alpha)A(\alpha)\right) \label{eq:LMI_2} 
\end{equation}  
have all positive definite coefficients.

Let $P(\alpha)$ be a homogeneous polynomial of degree $d_p$ which can be represented as 
\begin{equation}
P(\alpha)=\sum_{\gamma \in W_{d_p} } P_{\langle{\gamma}\rangle} \alpha^{{\gamma}}, \label{eq:P_alpha} 
\end{equation}
where the coefficients $P_{\langle \gamma \rangle} \in \mathbb{S}^n$ and where we recall that $W_{d_p}:=\left\lbrace \gamma \in \N^l\,:\, \sum_{i=1}^l \gamma_i = d_p \right\rbrace$ is the set of the exponents of all $l$-variate monomials of degree $d_p$. Since $A(\alpha)$ is a homogeneous polynomial of degree $d_a$, we can write it as 
\begin{equation}
A(\alpha)=\sum_{\gamma \in W_{d_a} } A_{\langle{\gamma}\rangle} \alpha^{{\gamma}}, \label{eq:A_alpha} 
\end{equation}
where the coefficients $A_{\langle \gamma \rangle} \in \mathbb{R}^{n \times n}$. By substituting~\eqref{eq:P_alpha} and~\eqref{eq:A_alpha} into~\eqref{eq:LMI_1} and~\eqref{eq:LMI_2} and defining $d_{pa}$ as the degree of $P(\alpha)A(\alpha)$, the conditions of Theorem~\ref{Polya's Theorem} can be represented in the form
\begin{equation}
\sum_{h \in W_{d_p}} \beta_{\langle{h}\rangle, \langle{\gamma}\rangle }P_{\langle{h}\rangle} \succ 0; \quad \gamma \in W_{d_p + d_1} \; \text{and} \label{eq:LMI_3} 
\end{equation}
\begin{equation} 
\hspace*{-0.06in} \sum_{h \in W_{d_p}} (H_{\langle{h}\rangle, \langle{\gamma}\rangle}^T P_{\langle{h}\rangle} + P_{\langle{h}\rangle} H_{\langle{h}\rangle, \langle{\gamma}\rangle}) \prec 0; \; \gamma \in W_{d_{pa} + d_2} \label{eq:LMI_4}.  
\end{equation}
Here $\beta_{\langle{h}\rangle,\langle{\gamma}\rangle}$ is defined to be the scalar coefficient which multiplies $P_{\langle{h}\rangle}$ in the $\langle{\gamma}\rangle$-th monomial of the homogeneous polynomial $\left( \sum_{i=1}^l \alpha_i \right)^{d_1} P(\alpha)$ using the lexicographical ordering. Likewise $H_{\langle{h}\rangle, \langle{\gamma}\rangle} \in \mathbb{R}^{n \times n}$ is the term which left or right multiplies $P_{\langle{h}\rangle}$ in the $\langle{\gamma}\rangle$-th monomial of $\left( \sum_{i=1}^l \alpha_i \right)^{d_2} \left( A^T(\alpha)P(\alpha)+P(\alpha)A(\alpha) \right)$ using the lexicographical ordering. For an intuitive explanation as to how these $\beta$ and $H$ terms are calculated, we consider a simple example. Precise formulae for these terms will follow the example. \vspace*{0.05in}

\textit{1) Example:} \textit{Calculating the $\beta$ and $H$ coefficients.}

Consider $
A(\alpha)=A_1 \alpha_1 + A_2 \alpha_2$ and $P(\alpha)=P_1 \alpha_1+P_2 \alpha_2$.
By expanding Equation~\eqref{eq:LMI_1} for $d_1=1$ we have  
\begin{equation}
(\alpha_1+\alpha_2)P(\alpha)= P_{1}\alpha_1^2+(P_{1}+P_{2})\alpha_1 \alpha_2 + P_{2} \alpha_2^2. 
\end{equation}
The $ \beta_{\langle{h}\rangle,\langle{\gamma}\rangle}$ terms are then extracted as 
\begin{equation}
\beta_{1,1}=1,\; \beta_{2,1}=0,\; \beta_{1,2}=1,\;
\beta_{2,2}=1, \;\beta_{1,3}=0, \;\beta_{2,3}=1. 
\end{equation}
Next, by expanding Equation~\eqref{eq:LMI_2} for $d_2=1$ we have 
\begin{align}
&(\alpha_1+\alpha_2) \left(A^T(\alpha)P(\alpha)+P(\alpha)A(\alpha) \right)= \left( A^T_{1}P_{1}+P_{1}A_{1} \right)\alpha_1^3 \nonumber \\
& + \left(  A^T_{1}P_{1}+P_{1} A_{1} +A^T_{2}P_{1}+P_{1}A_{2} +A^T_{1}P_{2}+P_{2}A_{1}\right) \alpha_1^2 \alpha_2 \nonumber \\
&+\left( A^T_{2}P_{1}+P_{1}A_{2}+A^T_{1}P_{2}+P_{2}A_{1}+A^T_{2}P_{2}+P_{2}A_{2} \right)\alpha_1 \alpha_2^2 \nonumber \\
& + \left(A^T_{2}P_{2}+P_{2}A_{2} \right)\alpha_2^3. 
\end{align}
The $ H_{\langle{h}\rangle,\langle{\gamma}\rangle} $ terms are then  extracted as
\begin{align}
&\hspace*{-0.1in} H_{1,1}=A_{1}, && \hspace*{-0.07in} H_{2,1}=\textbf{0},  && \hspace*{-0.4in} H_{1,2}=A_{1}+A_{2},    && \hspace*{-0.07in}  H_{2,2}=A_{1}, \nonumber \\ 
& \hspace*{-0.1in} H_{1,3}=A_{2}, &&  \hspace*{-0.07in} H_{2,3}=A_{1}+A_{2}, && \hspace*{-0.07in} H_{1,4}=\textbf{0},   && \hspace*{-0.07in} H_{2,4}=A_{2}.
\end{align}

\textit{2) General Formula:}
The $\{ \beta_{\langle{h}\rangle,\langle{\gamma}\rangle} \}$ can be formally defined recursively as follows. Let the initial values for $ \beta_{\langle{h}\rangle,\langle{\gamma}\rangle}$ be defined as 
\begin{equation}
\beta^{(0)}_{\langle{h}\rangle,\langle{\gamma}\rangle} = \begin{cases}1& \text{if } h=\gamma\\ 0 & \text{otherwise} \end{cases}\qquad \text{for} \; \gamma \in W_{d_p} \; \text{and} \; h \in W_{d_p}. \label{eq:beta_init} 
\end{equation}
Then, iterating for $i=1,\ldots d_1$, we let 
\begin{equation}
 \beta^{(i)}_{\langle{h}\rangle,\langle{\gamma}\rangle}=\sum_{\lambda \in W_1} \beta^{(i-1)}_{\langle{h}\rangle,\langle{\gamma-\lambda}\rangle} \qquad \text{for} \; \gamma \in W_{d_p + i} \; \text{and} \; h \in W_{d_p}. \label{eq:beta}
\end{equation}
Finally, we set $\{\beta_{\langle{h}\rangle,\langle{\gamma}\rangle}\} =  \{\beta^{d_1}_{\langle{h}\rangle,\langle{\gamma}\rangle}\}$. To obtain $\lbrace{ H_{\langle{h}\rangle, \langle{\gamma}\rangle}} \rbrace$, set the initial values as 
\begin{equation}
H^{(0)}_{\langle{h}\rangle,\langle{\gamma}\rangle} = \sum_{\lambda \in W_{d_a}: \lambda + h = \gamma} A_{\langle{\lambda} \rangle}
 \;\; \text{for} \; \gamma \in W_{d_p+d_a}\; \text{and} \; h \in W_{d_p}.  \label{eq:H_init} 
\end{equation}
Then, iterating for $i=1,\ldots d_2$, we let 
\begin{equation}
H^{(i)}_{\langle{h}\rangle,\langle{\gamma}\rangle}=\sum_{\lambda \in W_1 } H^{(i-1)}_{\langle{h}\rangle,\langle{\gamma-\lambda}\rangle}  \quad \text{for} \; \gamma \in W_{d_{pa} + i}\; \text{and} \;h \in W_{d_p}.  \label{eq:H} 
\end{equation}
Finally, set $\{H_{\langle{h}\rangle,\langle{\gamma}\rangle}\} = \{H^{d_2}_{\langle{h}\rangle,\langle{\gamma}\rangle}\}$.

\begin{figure*}[htbp]
\begin{minipage}[b]{0.5\linewidth} \hspace*{0.1in}
\includegraphics[scale=0.36]{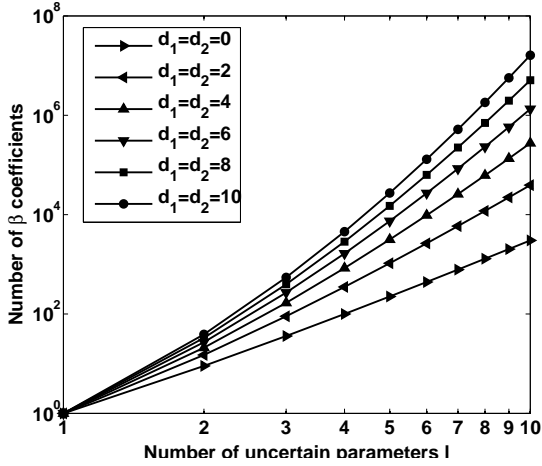}
\end{minipage} \hspace*{-0.1in} 
\begin{minipage}[b]{0.5\linewidth}  
 \includegraphics[scale=0.4]{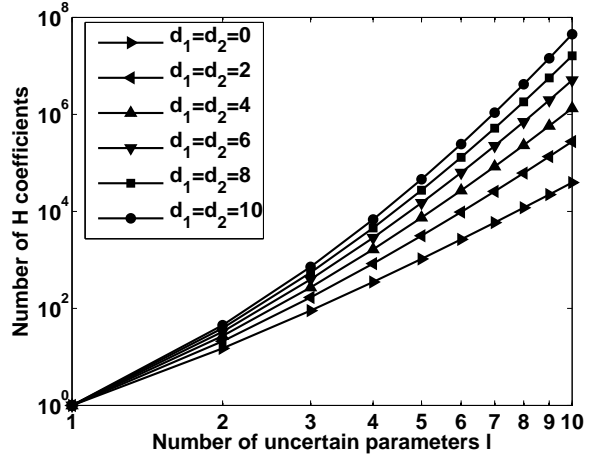}
\end{minipage} 
\caption{Number of $ \beta_{\langle{h}\rangle,\langle{\gamma}\rangle}$ and $ H_{\langle{h}\rangle,\langle{\gamma}\rangle}$ coefficients vs. the number of uncertain parameters for different Polya's exponents and for $d_p=d_a=2$}
\label{fig:No_coeff} 
\end{figure*}

For the case of large-scale systems, computing and storing $\lbrace \beta_{\langle{h}\rangle,\langle{\gamma}\rangle} \rbrace$ and $ \lbrace {H_{\langle{h}\rangle,\langle{\gamma}\rangle}} \rbrace$ is a significant challenge due to the number of these coefficients. Specifically, the number of terms increases with $l$ (number of uncertain parameters in system~\eqref{system}), $d_{p}$ (degree of $P(\alpha)$), $d_{pa}$ (degree of $P(\alpha)A(\alpha)$) and $d_1, d_2$ (Polya's exponents) as follows. \vspace*{0.1in}

\textit{3) Number of $ \beta_{\langle{h}\rangle,\langle{\gamma}\rangle}$ coefficients:}
For given $l,d_p$ and $d_1$, since $h \in W_{d_p}$ and $\gamma \in W_{d_p+d_1}$, the number of $ \beta_{\langle{h}\rangle,\langle{\gamma}\rangle}$ coefficients is the product of $L_0:= \text{card}(W_{d_p})$ and $L:=\text{card}(W_{d_p+d_1})$. Recall that card$(W_{d_p})$ is the number of all $l$-variate monomials of degree $d_p$ and can be calculated using~\eqref{eq:f} as follows. 
\begin{equation}
L_0 = f(l,d_p) =
\begin{cases}
\hspace*{0.2in}	0 &\text{for} \;\; l=0\\
\dbinom{d_p+l-1}{l-1}= \dfrac{(d_p+l-1)!}{d_p!(l-1)!} & \text{for} \;\; l > 0. \label{eq:L0} 
\end{cases}
\end{equation}
Likewise, card$(W_{d_p+d_1})$, i.e., the number of all $l-$variate monomials of degree $d_p+d_1$ is calculated using~\eqref{eq:f} as follows. 
\begin{align}
& \hspace*{-0.1in} L = f(l,d_p+d_1) = \nonumber \\
& \hspace*{-0.15in} \begin{cases}
\hspace*{0.2in}	0 &\text{for} \;\; l=0\\
\dbinom{d_p+d_1+l-1}{l-1}= \dfrac{(d_p+d_1+l-1)!}{(d_p+d_1)!(l-1)!} & \text{for} \;\; l > 0. \label{eq:L}
\end{cases}
\end{align}
The number of $ \beta_{\langle{h}\rangle,\langle{\gamma}\rangle}$ coefficients is $L_0 \cdot L$. \vspace*{0.1in}

\textit{4) Number of $ H_{\langle{h}\rangle,\langle{\gamma}\rangle} $ coefficients:}
For given $l,d_p, d_a$ and $d_2$, since $h\in W_{d_p}$ and $\gamma \in W_{d_{pa}+d_2}$, the number of $ H_{\langle{h}\rangle,\langle{\gamma}\rangle}$ coefficients is the product of $L_0:= \text{card}(W_{d_p})$ and $M:=\text{card}(W_{d_{pa}+d_2})$. By using~\eqref{eq:f}, we have 
\begin{align}
&M = f(l,d_{pa}+d_2) = \nonumber \\
&\begin{cases}
\hspace*{0.2in}	0 &\text{for} \;\; l=0\\
\dbinom{d_{pa}+d_2+l-1}{l-1}= \dfrac{(d_{pa}+d_2+l-1)!}{(d_{pa}+d_2)!(l-1)!} & \text{for} \;\; l > 0. \label{eq:M}
\end{cases}
\end{align}
The number of $ H_{\langle{h}\rangle,\langle{\gamma}\rangle}$ coefficients is $L_0 \cdot M$.

\begin{figure}[ht]
\centering
\includegraphics[scale=0.35]{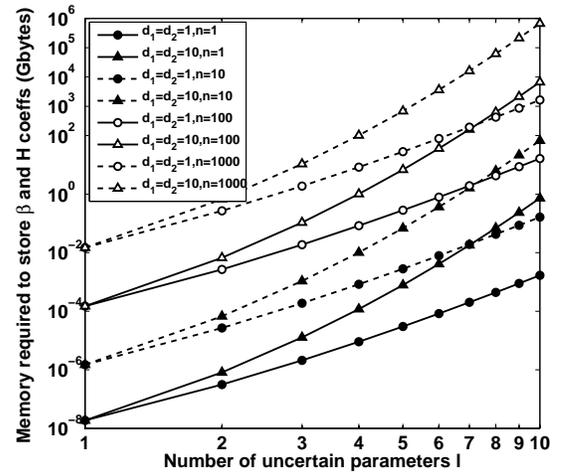} \vspace*{-0.15in}
\caption{Memory required to store $ \beta$ and $ H$ coefficients vs. number of uncertain parameters, for different $d_1,d_2$ and $d_p=d_a=2$}
\label{fig:memory_beta_h} 
\end{figure}

The number of $ \beta_{\langle{h}\rangle,\langle{\gamma}\rangle}$ and $ H_{\langle{h}\rangle,\langle{\gamma}\rangle}$ coefficients and the required memory to store these coefficients are shown in Figs.~\ref{fig:No_coeff} and~\ref{fig:memory_beta_h} in terms of the number of uncertain parameters $l$ and for different Polya's exponents. In all cases $d_p=d_a=2$.

It is observed from Fig.~\ref{fig:memory_beta_h} that, even for small $d_p$ and $d_a$, the required memory is in the Terabyte range.
In~\cite{peet_acc}, we proposed a decentralized computing approach to the
calculation of $\{ \beta_{\langle{h}\rangle,\langle{\gamma}\rangle} \}$ on large cluster computers. In the present work, we extend this method to the calculation of $ \{ H_{\langle{h}\rangle,\langle{\gamma}\rangle} \}$
and the SDP elements which will be discussed in the following section. We express the LMIs associated with conditions~\eqref{eq:LMI_3} and~\eqref{eq:LMI_4} as an SDP in both primal and dual forms. We also discuss the structure of the primal and dual SDP variables and the constraints.

\subsection {SDP Problem Elements}
\label{sec:structure}
A semi-definite programming problem can be stated either in primal or dual format.
Given $C \in \mathbb{S}^{m} $,  $a \in \mathbb{R}^K$ and $B_i \in \mathbb{S}^m$, the \textbf{primal problem} is of the form 
\begin{equation}
 \max_X \quad tr(CX) \nonumber 
\end{equation}
\begin{equation}
 \hspace*{0.1in} \text{subject to} \quad  a-B(X) = 0 \nonumber 
\end{equation}
\begin{equation}
\hspace*{0.35in} X\succeq 0,  \label{eq:primal} 
\end{equation}
where the linear operator $B: \mathbb{S}^{m} \rightarrow \mathbb{R}^K$ is defined as  
\begin{equation}
\hspace*{-0.1in}
B(X)= \left[
\begin{array}{ccc}
tr(B_1X) \quad tr(B_2 X) \quad \cdots \quad tr(B_K X) \end{array} \right]^T.  
\end{equation}
$X \in \mathbb{S}^{m}$ is the primal variable. Given a primal SDP, the associated \textbf{dual problem} is 
\begin{equation}
\hspace*{-0.2in}\min_{y,Z} \quad a^Ty \nonumber 
\end{equation}
\begin{equation}
 \hspace*{0.2in}\text{subject to} \quad  B^T(y)-C = Z \nonumber 
\end{equation}
\begin{equation}
 \hspace*{1in}Z\succeq 0 \;,\; y \in \mathbb{R}^K,  \label{eq:dual} 
\end{equation}
where $B^T: \mathbb{R}^K \rightarrow \mathbb{S}^{m}$ is the transpose operator and is given by 
\begin{equation}
B^T(y)=\sum_{i=1}^K y_i B_i  \label{eq:AT}
\end{equation}
and where $y \in \mathbb{R}^K$ and $Z \in \mathbb{S}^{m}$ are the dual variables.
The elements $C$, $B_i$ and $a$ of the SDP problem associated with the LMIs in~\eqref{eq:LMI_3} and~\eqref{eq:LMI_4} are defined as follows.
We define the element $C$ as 
\begin{equation}
C := \text{diag}(C_1, \cdots C_L, C_{L+1}, \cdots C_{L+M}), \label{eq:C}
\end{equation}
where 
\begin{equation}
C_i :=
\begin{cases}
	\delta I_n \cdot \left( \sum_{h \in W_{d_p}} \beta_{\langle h \rangle,i} \, \frac{d_p!}{h_1! \, \cdots \, h_l!} \right), & \hspace*{0.35in} 1 \le i \le L\\
	0_n, & \hspace*{-0.25in}   L+1 \le  i  \le L+M,
\end{cases}  \label{eq:Cj}
\end{equation}
where recall that $L= \text{card}(W_{d_p+d_1})$ is the number of monomials in $\left(\sum_{i=1}^l \alpha_i \right)^{d_1} P(\alpha)$, $M=\text{card}(W_{d_{pa}+d_2})$ is the number of monomials in $\left(\sum_{i=1}^l \alpha_i \right)^{d_2} P(\alpha)A(\alpha)$, where  $n$ is the dimension of system~\eqref{system}, $l$ is the number of uncertain parameters and $\delta$ is a small positive parameter.

For $i=1, \cdots, K$, define $B_i$ elements as
\begin{equation}
B_i := \text{diag}(B_{i,1}, \cdots B_{i,L}, B_{i,L+1}, \cdots B_{i,L+M}), \label{eq:Ai} 
\end{equation}
where $K$ is the number of dual variables in (30) and is equal to the product of  the number of upper-triangular elements in each $P_\gamma \in \mathbb{S}^n$ (the coefficients in $P(\alpha)$) and the number of coefficients in $P(\alpha)$ (i.e. the cardinality of $W_{d_p}$). Since there are $f(l,d_p)=\dbinom{d_p+l-1}{l-1}$ coefficients in $P(\alpha)$ and each coefficient has $\tilde{N}:=\frac{1}{2}n(n+1)$ upper-triangular elements, we find 
\begin{equation}
K= \frac{(d_p + l-1)!}{d_p! (l-1)!}\tilde{N}. \label{eq:K} 
\end{equation}
To define the $B_{i,j}$ blocks, first we define the function $V_{\langle h \rangle}: \mathbb{Z}^K \rightarrow \mathbb{Z}^{n \times n}$, 
\begin{equation}
V_{\langle h \rangle}(x) := \sum_{j=1}^{\tilde{N}} E_j \; x_{j+\tilde{N}(\langle h \rangle-1)} \quad \text{for all} \quad h \in W_{d_p}, 
\end{equation}
which maps each variable to a basis matrix $E_j$, where recall that $E_j$ is the basis for $ \mathbb{S}^n $. Note that a different choice of basis would require a different function $V_{\langle h \rangle}$. Then for $i=1, \cdots, K$,  
\begin{align}
&B_{i,j}:= \nonumber \\
&\begin{cases}
\sum_{h \in W_{d_p}} \beta_{\langle h \rangle,j} V_{\langle h \rangle}(e_i), \hspace*{1in} 1 \le j \le L \;  (I)\\
-\sum_{h \in W_{d_p}} \Big( H_{{\langle h \rangle} ,j-L}^T V_{\langle h \rangle}(e_i)+  \\
\hspace*{0.6in} V_{\langle h \rangle}(e_i) H_{{\langle h \rangle},j-L} \Big), \;\;  L+1 \le j \le L+M. \; (II)
\end{cases} \label{eq:Aij} 
\end{align}
Finally, to complete the SDP problem associated with Polya's algorithm set
\begin{equation}
a=\vec{1} \in \mathbb{R}^K. \label{eq:a} 
\end{equation}

\subsection {Parallel Set-up Algorithm}
\label{sec:set-up parallel}

In this section, we propose a decentralized, iterative algorithm for calculating the terms $\{ \beta_{\langle h \rangle, \langle \gamma \rangle} \}$, $\{H_{\langle h \rangle, \langle \gamma \rangle}\}$, $C$ and $B_{i}$ as defined in~\eqref{eq:beta},~\eqref{eq:H},~\eqref{eq:C} and~\eqref{eq:Ai}. The algorithm has been implemented in C++, using MPI (Message Passing Interface) and is available at: www.sites.google.com/a/asu.edu/kamyar/software. We present an abridged description of this algorithm in Algorithm~\ref{alg:setup}, wherein $N$ is the number of available processors.

\begin{algorithm}
\vspace*{0.05in}
\textbf{\textit{Inputs}:} $d_p$: degree of $P(\alpha)$, $d_a$: degree of $A(\alpha)$, $n$: number of states, $l$: number of uncertain parameters, $d_1, d_2$: number of Polya's iterations, Coefficients of $A(\alpha)$.

\textbf{\textit{Initialization}}:
Set $\hat{d}_1=\hat{d}_2=0$ and $d_{pa}=d_p+d_a$. Calculate $L_0$ as the number of monomials in $P(\alpha)$ using~\eqref{eq:L0} and $M$ as the number of monomials in $P(\alpha)A(\alpha)$ using \eqref{eq:M}. Set $L=L_0$. Calculate $ L'=\mathtt{floor}(\frac{L}{N})$ and $M'=\mathtt{floor}(\frac{M}{N})$ as the number of monomials in $P(\alpha)$ and $P(\alpha)A(\alpha)$ assigned to each processor.\\
 \For{$i=1, \cdots,N$, processor $i$}{
 Initialize $\beta_{k,j}$ for $j=(i-1)L'+1, \cdots, iL'$ and $k=1, \cdots L_0$ using~\eqref{eq:beta_init}.\\
 Initialize $H_{k,m}$ for $m=(i-1)M'+1, \cdots, iM'$ and $k=1, \cdots L_0$ using~\eqref{eq:H_init}.
 } 

 \textbf{\textit{Calculating $\beta$ and $H$ coefficients:}} \\
 \While{$\hat{d}_1 \leq d_1 $ or $\hat{d}_2 \leq d_2$}{  
  \If{$\hat{d}_1 \leq d_1 $}{ 
  \For{$i=1, \cdots,N$, processor $i$}{
  Set $d_p=d_p+1$ and $\hat{d}_1=\hat{d}_1+1$. Update $L$ using~\eqref{eq:L}. Update $L'=\mathtt{floor}(\frac{L}{N})$.\\
 Calculate $\beta_{k,j}$ for $j=(i-1)L'+1, \cdots, iL'$ and $k=1, \cdots L_0$ using~\eqref{eq:beta}.}
 } 
 \If{$\hat{d}_2 \leq d_2 $}{
  \For{$i=1, \cdots,N$, processor $i$}{
 \hspace*{-0.1in}  Set $d_{pa}=d_{pa} \hspace*{-0.05in} +  \hspace*{-0.05in}1$ and $\hat{d}_2=\hat{d}_2  \hspace*{-0.05in}+  \hspace*{-0.05in}1$. Update $M$ using~\eqref{eq:M}. Update $M'=\mathtt{floor}(\frac{M}{N})$.\\
Calculate $H_{k,m}$ for $m=(i-1)M'+1, \cdots, iM'$ and $k=1, \cdots L_0$  using~\eqref{eq:H}.} 
 } 
 } 

 \textbf{\textit{Calculating the SDP elements}:} \\
 \For{$i=1, \cdots,N$, processor $i$}{
  Calculate the number of dual variables $K$ using~\eqref{eq:K}. Set $T' = \mathtt{floor}(\frac{L+M}{N})$.\\
 Calculate the blocks of the SDP element $C$ as 
\[
 \begin{cases}
 C_j \; \text{using~\eqref{eq:Cj}} & \text{for} \; j = (i-1)L'+1, \cdots, iL' \\
 C_j = 0_n & \hspace*{-0.2in} \text{for} \; j=L+(i-1)M'+1, \cdots, L+iM'
 \end{cases} 
 \]
  Set the sub-blocks of the SDP element $C$ as 
  \begin{equation}
 \overline{\textbf{C}}_i = \text{diag}\left(C_k|_{k=(i-1)T'+1}^{iT'}\right). \label{eq:Cbar}  
  \end{equation}
 \For{$j=1, \cdots,K$}{
 Calculate the blocks of the SDP elements $B_j$ as 
\[
 \begin{cases}
 B_{j,k}  \; \text{using~\eqref{eq:Aij}-\textit{I}} & \text{for} \; k=(i-1)L'+1, \cdots, iL' \\
 B_{j,k} \; \text{using~\eqref{eq:Aij}-\textit{II}} & \text{for} \; k=L+(i-1)M'+1, \\ & \qquad \qquad \qquad \quad \cdots, L+iM'
 \end{cases}
 \]
  Set the sub-blocks of the SDP element $B_j$ as 
  \begin{equation}
\overline{\textbf{B}}_{j,i}= \text{diag} \left(B_{j,k}|_{k=(i-1)T'+1}^{iT'} \right)  . \label{eq:Bbar}   
  \end{equation}
 }
 }

 \textbf{\textit{Outputs}:}  Sub-blocks  $\overline{\textbf{C}}_i$ and $\overline{\textbf{B}}_{j,i}$ of the SDP elements for $i=1,\cdots,N$ and $j=1,\cdots,K$.
\caption{The parallel set-up algorithm} 
\label{alg:setup}
\end{algorithm}

Note that we have only addressed the problem of robust stability analysis, using the polynomial inequality
\[
P(\alpha) \succ 0,  A^T(\alpha)P(\alpha) + P(\alpha) A(\alpha) \prec 0
\]
 for $ \alpha \in \Delta_l$. However, we can generalize the decentralized set-up algorithm to consider a more general class of feasibility problems, i.e., 
\begin{equation}
\sum_{i=1}^{\hat{N}} \left( \tilde{A}_i(\alpha)\tilde{X}(\alpha)\tilde{B}_i(\alpha) + \tilde{B}^T_i(\alpha)\tilde{X}(\alpha)\tilde{A}^T_i(\alpha) + R_i(\alpha) \right) \prec 0  \label{eq:LMI_general} 
\end{equation}
for $\alpha \in \Delta_l$. 
One motivation behind the development of such generalized set-up algorithm is that the parameter-dependent versions of the LMIs associated with $H_2$ and $H_\infty$ synthesis problems in~\cite{Gahinet94alinear,dullerud2000course} can be formulated in the form of~\eqref{eq:LMI_general}. 

\subsection {Set-up algorithm: Complexity Analysis}

Since checking the positive definiteness of all representatives of a square matrix with parameters on proper real intervals is intractable~\cite{nemirovskii1993several}, the question of feasibility of~\eqref{eq:Lyap_LMI} is also intractable. To solve the problem of inherent intractability we establish a trade off between accuracy and complexity. In fact, we develop a sequence of decentralized polynomial-time algorithms whose solutions converge to the exact solution of the NP-hard problem. In other words, the translation of a polynomial optimization problem to an LMI problem is the main source of complexity. This high complexity is unavoidable and, in fact, is the reason we seek parallel algorithms.


Algorithm~\ref{alg:setup} distributes the computation and storage of $\{ \beta_{\langle h \rangle, \langle \gamma \rangle} \}$ and $\{H_{\langle h \rangle, \langle \gamma \rangle}\}$ among the processors and their dedicated memories, respectively. In an ideal case, where the number of available processors is sufficiently large (equal to the number of monomials in $P(\alpha) A(\alpha)$, i.e. $M$) only one monomial ($L_0$ of $\beta_{\langle h \rangle, \langle \gamma \rangle}$ and $L_0$ of $H_{\langle h \rangle, \langle \gamma \rangle}$) are assigned to each processor. \vspace*{0.05in}

\textit{1) Computational complexity analysis:}
The most computationally expensive part of the set-up algorithm is the calculation of the $B_{i,j}$ blocks in~\eqref{eq:Aij}. Considering that the cost of matrix-matrix multiplication is $\sim n^3$, the cost of calculating each $B_{i,j}$ block is
$
\sim \text{card}(W_{d_p}) \cdot n^3.
$
According to~\eqref{eq:Ai} and~\eqref{eq:Aij}, the total number of $B_{i,j}$ blocks is $K(L+M)$. Hence, as per Algorithm 1, each processor processes $K \left( \mathtt{floor}(\frac{L}{N} ) + \mathtt{floor}(\frac{M}{N}) \right)$ of the $B_{i,j}$ blocks, where $N$ is the number of available processors. Thus the per processor computational cost of calculating the $B_{i,j}$ at each Polya's iteration is 
\begin{equation}
\sim \text{card}(W_{d_p}) \cdot n^3 \cdot  K \left( \mathtt{floor}\left(\frac{L}{N}\right) + \mathtt{floor}\left(\frac{M}{N} \right) \right). 
\end{equation}
By substituting for $K$ from~\eqref{eq:K}, card$(W_{d_p})$ from~\eqref{eq:L0}, $L$ from~\eqref{eq:L} and $M$ from~\eqref{eq:M}, the per processor computation cost at each iteration is
\begin{align}
\begin{split}
\sim &\left( \frac{(d_p + l-1)!}{d_p! (l-1)!} \right)^2 \hspace*{-0.05in}  \frac{n^4}{2(n+1)} \left( \mathtt{floor}  \left( \frac{ \dfrac{(d_p+d_1+l-1)!}{(d_p+d_1)!(l-1)!}}{N} \right) \right. \\
  & \left. +  \mathtt{floor} \left( \frac {\dfrac{(d_{pa}+d_2+l-1)!}{(d_{pa}+d_2)!(l-1)!}} {N} \right) \vphantom{\frac12}   \vphantom{\frac12}\right)
\end{split}
\end{align}
assuming that $l > 0$ and $N \leq M$. For example, for the case of large-scale systems (large $n$  and $l$), the computation cost per processor at each iteration is $\sim (l^{2d_p+d_1}+l^{2d_p+d_a+d_2}) n^5$ having $N=L_0$ processors, $\sim (l^{2d_p+d_1}+l^{2d_p+d_a+d_2}) n^5$ having $N=L$ processors and $\sim l^{2d_p+d_a+d_2-d_1}n^5$ having $N=M$ processors. Thus for the case where $d_p \geq 3$, the number of operations grows more slowly in $n$ than in $l$. \vspace*{0.05in}

\begin{figure}[t]
\centering
\subfigure[]{
  \includegraphics[scale=0.27]{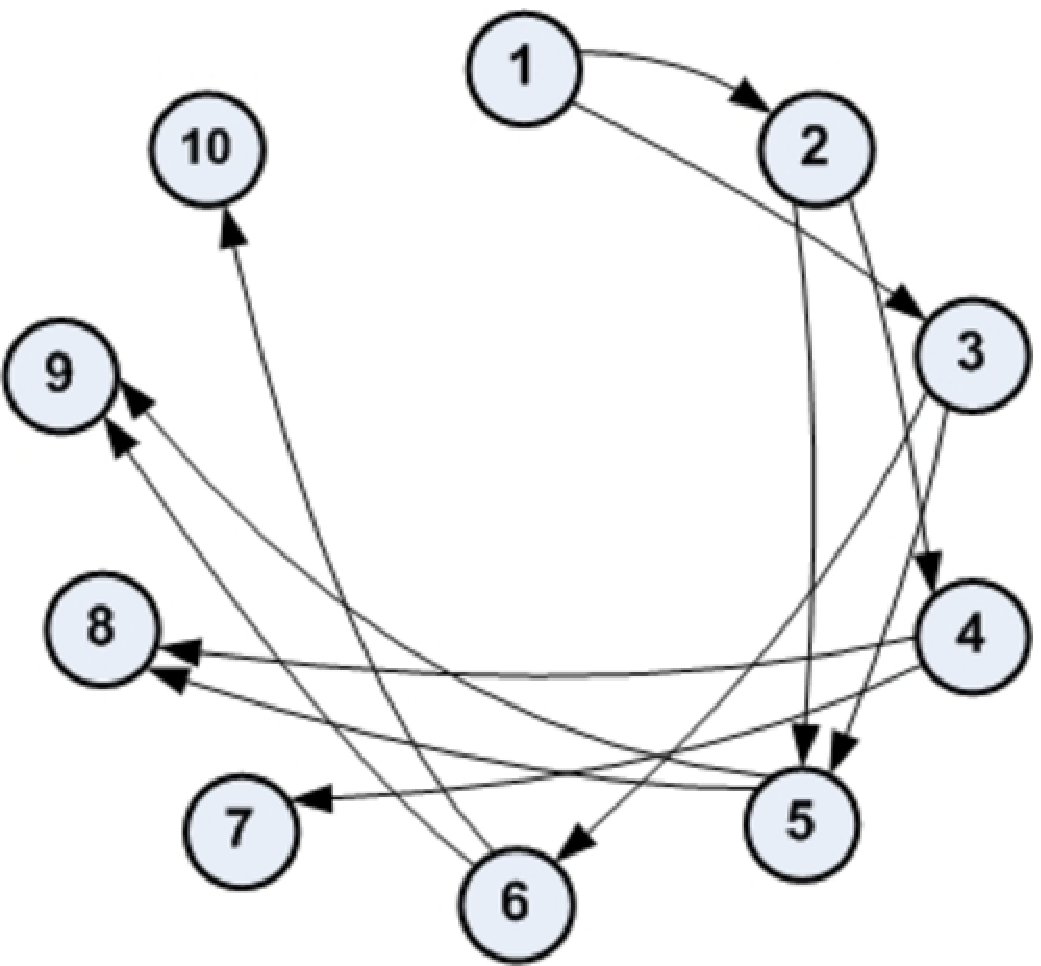}
	    \label{fig:graph1}
}  \hspace*{0.1in}
\subfigure[]{
\includegraphics[scale=0.35]{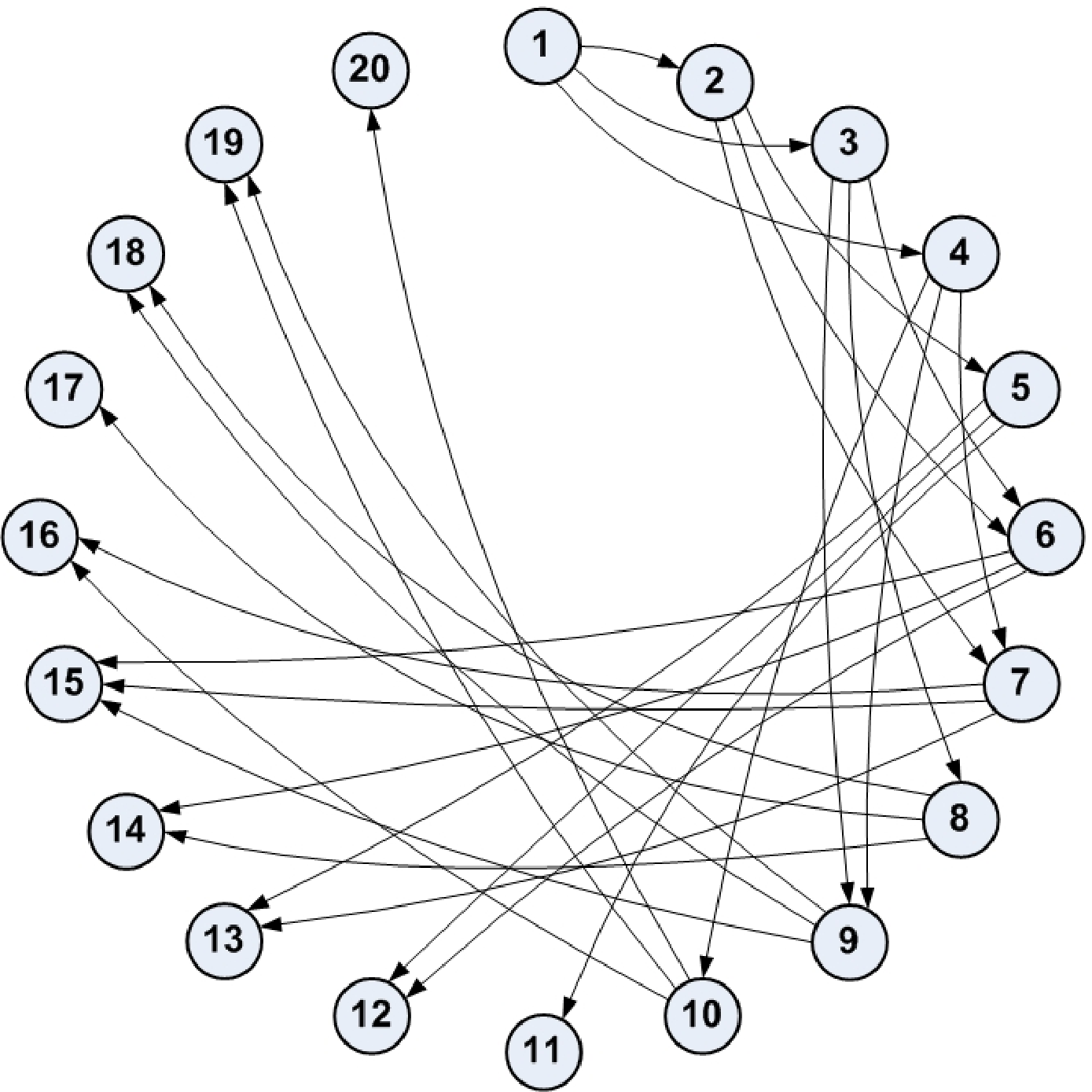}
	    \label{fig:graph2} 
}
\caption{Graph representation of the network communication of the set-up algorithm. (a) Communication directed graph for the case $\alpha \in \Delta_3,d_p=2$. (b) Communication directed graph for the case $\alpha \in \Delta_4,d_p=2$. } 
\end{figure}

\textit{2) Communication complexity analysis:}
Communication between processors can be modeled by a directed graph $G(V,E)$, where the set of nodes $V=\{ 1, \cdots, N\}$ is the set of indices of the available processors and the set of edges $E=\{ (i,j):i,j \in V \}$ is the set of all pairs of processors that communicate with each other. For every directed graph we can define an adjacency matrix $T_G$. If processor $i$ communicates with processor $j$, then $[T_G]_{i,j} = 1$, otherwise $[T_G]_{i,j} = 0$.  In this section, we only define the adjacency matrix for the part of the algorithm that performs Polya's iterations on $P(\alpha)$. For Polya's iterations on $P(\alpha)A(\alpha)$, the adjacency matrix can be defined in a similar manner.
For simplicity, we assume that at each iteration, the number of available processors is equal to the number of monomials in $(\sum_{i=1}^l \alpha_i)^{d_1}P(\alpha)$. Using~\eqref{eq:L}, let us define $r_{d_1}$ and $r_{d_1+1}$ as the numbers of monomials in $(\sum_{i=1}^l \alpha_i)^{d_1}P(\alpha)$ and $(\sum_{i=1}^l \alpha_i)^{d_1+1}P(\alpha)$.
For $I=1, \cdots, r_{d_1}$, define
\begin{align*}
\mathcal{E}_I := \{ \text{lex. indices of monomials in }\; & \left( \sum_{i=1}^l \alpha_i \right) \alpha^{\gamma}: \\
 & \hspace*{-0.2in} \gamma \in W_{d_p+d_1} \; \text{and} \; \langle \gamma \rangle = I \}. 
\end{align*}
Then for $i=1, \cdots, r_{d_1+1}$ and $j=1, \cdots, r_{d_1+1}$,
\[
[T_G]_{i,j} :=
\begin{cases}
1 & \text{if} \quad i \leq r_{d_1} \; \text{and} \; j \in \mathcal{E}_i \; \text{and}   \; i \neq j \\
0 & \text{otherwise}.
\end{cases}
\]
Note that this definition implies that the communication graph of the set-up algorithm changes at every iteration.
To help visualize the graph, the adjacency matrix for the case where $\alpha \in \Delta_2$ is
\begin{small}
\[
T_G \hspace*{-0.025in} :=  \hspace*{-0.05in} \left[ \hspace*{-0.025in} 
\begin{array}{cccccc|ccc}
0 & 1 & 0 & \cdots &  0 &                          0 & 0 & \cdots & 0 \\
0 & 0 & 1 & 0 & \cdots  &                          0 & \\
\vdots & \vdots  & \ddots & \ddots & \ddots & \vdots & \vdots & \ddots & \vdots \\
\vdots & \vdots &  & \ddots & \ddots &             0 &  \\
0 & 0 & \cdots & \cdots & 0 &                      1 & 0 & \cdots & 0 \\
\hline
0 &  &  & \cdots & &                               0 & 0 & \cdots & 0 \\
\vdots &  &  & \ddots & &   \vdots   &  \vdots & \ddots &  \vdots \\
0 &  &  & \cdots & &                               0 &  0 & \cdots & 0\\
\end{array} \hspace*{-0.025in}  \right] \hspace*{-0.025in} \in \hspace*{-0.025in} \mathbb{R}^{r_{d_1+1} \times r_{d_1+1}},
\]
\end{small}
 where the nonzero sub-block of $T_G$ lies in $\mathbb{R}^{r_{d_1} \times r_{d_1}}$.
We can also illustrate the communication graphs for the cases $\alpha \in \Delta_3$ and $\alpha \in \Delta_4$ with $d_p=2$ as seen in Fig.~\ref{fig:graph1} and~\ref{fig:graph2}.

For a given algorithm, the communication complexity is defined as the sum of the size of all communicated messages. For simplicity, let us consider the worst case scenario, where each processor is assigned more than one monomial and sends all of its assigned $\beta_{\langle h \rangle, \langle \gamma \rangle} $ and $H_{\langle h \rangle, \langle \gamma \rangle} $ coefficients to other processors. In this case, the algorithm assigns $\mathtt {floor}(\frac{L}{N})\cdot \text{card}(W_{d_p})$ of the $\beta_{\langle h \rangle, \langle \gamma \rangle}$ coefficients, each of size 1, and $\left( \mathtt {floor}(\frac{L}{N}) + \mathtt {floor}(\frac{M}{N}) \right) \cdot \text{card}(W_{d_p})$ of the $H_{\langle h \rangle, \langle \gamma \rangle}$ coefficients, each of size $n^2$, to each processor. Thus the communication complexity of the algorithm per processor and per iteration is  
\begin{equation}
\text{card}(W_{d_p}) \left( \mathtt {floor}\left(\frac{L}{N}\right) + \mathtt {floor}\left(\frac{M}{N}\right) n^2 \right). \label{eq:setup_comm}
\end{equation}
This indicates that increasing the number of processors (up to $M$) actually leads to less communication overhead per processor and improves the scalability of the algorithm. By substituting for card$(W_{d_p})$ from~\eqref{eq:L0}, $L$ from~\eqref{eq:L} and $M$ from~\eqref{eq:M} and considering large $l$ and $n$, the communication complexity per processor at each Polya's iteration is $\sim  l^{d_{pa}+d_2} n^2 $ having $N=L_0$ processors, $\sim l^{d_{pa}+d_2-d_1} n^2$ having $N=L$ processors and $\sim l^{d_p} n^2$ having $N=M$ processors.

\section{PARALLEL SDP SOLVER}
\label{sec:SDPSOLVER}
In this section, we describe the steps of our primal-dual interior-point algorithm and show how, for the LMIs in~\eqref{eq:LMI_3} and~\eqref{eq:LMI_4}, these steps can be distributed in a distributed-computing, distributed-memory environment. 

\subsection {Interior-point methods}
\label{sec:Interior-point} 
Interior-point methods define a popular class of algorithms for solving linear and semi-definite programming problems.
The most widely accepted interior-point algorithms are dual scaling~\cite{dual_scaling,primal}, primal-dual~\cite{primal_dual_1,primal_dual_2,primal_dual_3} and cutting-plane/spectral bundle~\cite{helmberg2000spectral, krishnan, mahdu}. In this paper, we use the central-path-following primal-dual algorithm described in~\cite{primal_dual_3} and~\cite{csdp}. Although we found it possible to use dual-scaling algorithms, we chose to pursue a primal-dual algorithm because, in general, primal-dual algorithms converge faster~\cite{primal_dual_3,primal} while still preserving the structure of the solution (see~\eqref{eq:S}) at each iteration. We prefer primal-dual to cutting plane/spectral bundle methods because, as we show in Section~\ref{sec:complexity}, the centralized part of our primal-dual algorithm consists of solving a symmetric system of linear equations (see~\eqref{eq:SAE1}), whereas for the cutting plane/spectral bundle algorithm, the centralized computation would consist of solving a constrained quadratic program (see~\cite{krishnan},~\cite{mahdu}) with number of variables equal to the size of the system of linear equations. Because centralized computation is the limiting factor in a parallel algorithm, and because solving symmetric linear equations is simpler than solving a quadratic programming problem, we chose the primal-dual approach.

The choice of a central path-following primal-dual algorithm as in~\cite{primal_dual_3} and~\cite{Alizadeh_method} was motivated by results in~\cite{Alizadeh_comparison} which demonstrated better convergence, accuracy and robustness over the other types of primal-dual algorithms. More specifically, we chose the approach in~\cite{primal_dual_3} over~\cite{Alizadeh_method} because unlike the Schur complement matrix (SCM) approach of the algorithm in~\cite{Alizadeh_method}, the SCM of~\cite{primal_dual_3} is symmetric and only the upper-triangular elements need to be sent/received by the processors. This leads to less communication overhead. The other reason for choosing~\cite{primal_dual_3} is that the symmetric SCM of the algorithm in~\cite{primal_dual_3} can be factorized using Cholesky factorization, whereas the non-symmetric SCM of~\cite{Alizadeh_method} must be factorized by LU factorization (LU factorization is roughly twice as expensive as Cholesky  factorization). Since factorization of SCM comprises the main portion of centralized computation in our algorithm, it is crucial for us to use computationally cheaper factorization methods to achieve better scalability.

In the primal-dual algorithm, both primal and dual problems are solved by iteratively calculating primal and dual step directions and step sizes, and applying these to the primal and dual variables. Let $X$ be the primal variable and $y$ and $Z$ be the dual variables. At each iteration, the variables are updated as
\begin{align}
X_{k+1} & = X_k + t_p \Delta X  \label{eq:X} \\
y_{k+1} & = y_k + t_d \Delta y  \label{eq:y} \\
Z_{k+1} & = Z_k + t_d \Delta Z  \label{eq:Z},
\end{align}
where $\Delta X$, $\Delta y$, and $\Delta Z$ are Newton's search direction and $t_p$ and $t_d$ are primal and dual step sizes. We choose the step sizes using a standard line-search between $0$ and $1$ with the constraint that $X_{k+1}$ and $Z_{k+1}$ remain positive semi-definite. We use a Newton's search direction given by  
\begin{align}
&\Delta X = \Delta \widehat{X} + \Delta \overline{X}  \label{eq:totalstepX} \\
&\Delta y = \Delta \widehat{y} + \Delta \overline{y}  \label{eq:totalstepy} \\
&\Delta Z = \Delta \widehat{Z} + \Delta \overline{Z}, \label{eq:totalstepZ}
\end{align}
where $ \Delta \widehat{X}$, $\Delta \widehat{y}$ and $\Delta \widehat{Z}$ are the predictor step directions and $\Delta \overline{X}$, $\Delta \overline{y}$, and $\Delta \overline{Z}$ are the corrector step directions. As per~\cite{primal_dual_3}, the predictor step directions are found as  
\begin{align}
&\Delta \widehat{y} = \Omega^{-1}\left(-a + B( Z^{-1} G X )\right)   \\
&\Delta \widehat{X} = -X + Z^{-1} G B^T( \Delta \widehat{y} )X      \label{eq:delX_hat} \\
&\Delta \widehat{Z} = B^T(y) - Z - C + B^T ( \Delta \widehat{y}),   \label{eq:delZ_hat}
\end{align}
where $C$ and the operators $B$ and $B^T$ are as defined in the previous section,  
\begin{equation}
 G = -B^T ( y ) + Z + C,  \label{eq:G} 
 \end{equation}
 and  
 \begin{equation}
 \Omega = \left[ B( Z^{-1} B^T( e_1 ) X ) \; \cdots \; B( Z^{-1} B^T ( e_k ) X ) \right].
\end{equation}
Recall that $ e_1, ...,e_k $ are the standard basis for $ \mathbb{R}^k$. Once we have the predictor step directions, we can calculate the corrector step directions as per~\cite{primal_dual_3}. Let $\mu=\dfrac{1}{3}tr(ZX)$. The corrector step directions are 
\begin{align}
& \Delta \overline{y} = \Omega^{-1}\left(  B( \mu Z^{-1} ) - B( Z^{-1} \Delta \widehat{Z} \Delta \widehat{X} ) \right)  \\
& \Delta \overline{X}   =   \mu Z^{-1}-Z^{-1} \Delta \widehat{Z} \Delta \widehat{X} - Z^{-1} \Delta \overline{Z} X \label{eq:delx_bar}  \\
& \Delta \overline{Z}   =   B^T( \Delta \overline {y}). \label{eq:delz_bar}   
\end{align}
The stopping criterion is $|a^T y - tr(CX)|\le \epsilon$. Information regarding the selection of starting points and convergence of different variants of interior-point primal-dual algorithm, including the algorithm we use in this paper are presented in~\cite{primal_dual_1},~\cite{primal_dual_2} and~\cite{primal_dual_3}.
\vspace*{-0.1in}

\subsection {Structure of SDP Variables} 
The key algorithmic insight of this paper which allows us to use the primal-dual approach presented in~\cite{primal_dual_3} is that by choosing an initial value for the primal variable with a certain block structure corresponding to the distributed structure of the processors, the algorithm will preserve this structure at every iteration. Specifically, we define the following structured block-diagonal subspace where each block corresponds to a single processor. 
\begin{align}
S_{l,m,n}& := \{ Y \subset \mathbb{R}^{(l+m)n \times (l+m)n} : \nonumber \\
& \hspace*{-0.25in} Y=   \text{diag}(Y_1,\cdots Y_l, Y_{l+1}, \cdots Y_{l+m}) \; \text{for} \; Y_i \in \mathbb{R}^{n \times n} \} \label{eq:S}
\end{align}
According to the following theorem, the subspace $S_{l,m,n}$ is invariant under Newton's iteration in the sense that when the algorithm in~\cite{primal_dual_3} is applied to the SDP problem defined by the polynomial optimization problem with initial value of the primal variable $X_0 \in S_{l,m,n}$, then the primal variable remains in the subspace at every Newton's iteration $X_k$.\\ \vspace*{-0.1in}

\begin{thm}\label{thm:thm3}
Consider the SDP problem defined in~\eqref{eq:primal} and~\eqref{eq:dual} with elements given by~\eqref{eq:C},~\eqref{eq:Ai} and~\eqref{eq:a}. Suppose $L$ and $M$ are the cardinalities of $W_{d_p+d_1}$ and $W_{d_{pa} +d_2}$. If~\eqref{eq:X},~\eqref{eq:y} and~\eqref{eq:Z} are initialized by 
\begin{equation}
X_0 \in S_{L,M,n}, \quad y_0 \in \mathbb{R}^{K}, \quad Z_0 \in S_{L,M,n}, 
\end{equation}
then for all $ k \in \mathbb{N} $, 
\begin{equation}
 X_k \in S_{L,M,n}, \quad Z_k \in S_{L,M,n}. 
\end{equation}
\end{thm}
\vspace*{0.05in}
\begin{proof}
We proceed by induction. First, suppose for some $k \in \mathbb{N}, $ 
\begin{equation}
 X_k \in S_{L,M,n} \quad \text{and}\quad Z_k \in S_{L,M,n}. \label{eq:assume} 
\end{equation}
We would like to show that this implies $X_{k+1}, Z_{k+1} \in S_{L,M,n}$. To see this, observe that according to~\eqref{eq:X} 
\begin{equation}
X_{k+1} = X_k + t_p \Delta X_k \quad \text{for all} \;\, k \in \mathbb{N}. 
\end{equation}
From~\eqref{eq:totalstepX}, $\Delta X_k$ can be written as  
\begin{equation}
\Delta X_k = \Delta \widehat{X}_k + \Delta \overline{X}_k \quad \text{for all} \;\, k \in \mathbb{N}.  \label{eq:DelX_k}
\end{equation}
To find the structure of $\Delta X_k$, we focus on the structures of $\Delta \widehat{X}_k$ and $\Delta \overline{X}_k$ individually. Using~\eqref{eq:delX_hat}, $\Delta \widehat{X}_k$ is 
\begin{equation}
\Delta \widehat{X}_k = -X_k + Z_k^{-1} G_k B^T( \Delta \widehat{y}_k )X_k \quad \text{for all} \;\, k \in \mathbb{N}. \label{eq:Dxhat_k} 
\end{equation}
where according to~\eqref{eq:G}, $G_k$ is
\begin{equation}
G_k = C-B^T ( y_k ) + Z_k \quad \text{for all} \;\, k \in \mathbb{N}.  \label{eq:G_k1}
\end{equation}
First we examine the structure of $G_k$. According to the definition of $C$ and $B_i$ in~\eqref{eq:C} and~\eqref{eq:Ai}, and the definition of  $B^T(y)$ in~\eqref{eq:AT}, we know that 
\begin{equation}
C \in S_{L,M,n}, \quad B^T : \R^K \mapsto S_{L,M,n}.  \label{eq:CA} 
\end{equation}
Since all the terms on the right hand side of~\eqref{eq:G_k1} are in $S_{L,M,n}$ and $S_{L,M,n}$ is a subspace, we conclude 
\begin{equation}
G_k \in S_{L,M,n}. \label{eq:G_k} \vspace*{-0.05in}
\end{equation}
Returning to~\eqref{eq:Dxhat_k}, using our assumption in~\eqref{eq:assume} and noting that the structure of the matrices in $S_{L,M,n}$ is also preserved through multiplication and inversion, we conclude 
\begin{equation}
\Delta \widehat{X}_k \in S_{L,M,n}. \label{eq:delxhat}
\end{equation}
Using~\eqref{eq:delx_bar}, the second term in~\eqref{eq:DelX_k} is 
\begin{equation}
\Delta \overline{X}_k = \mu Z_k^{-1}-Z_k^{-1} \Delta \widehat{Z}_k \Delta \widehat{X}_k - Z_k^{-1} \Delta \overline{Z}_k X_k \quad \text{for all} \;\, k \in \mathbb{N}. \label{Delxbar} 
\end{equation}
To determine the structure of $\Delta \overline{X}_k$, first we investigate the structure of $\Delta \widehat{Z}_k$ and $\Delta \overline{Z}_k$. According to~\eqref{eq:delZ_hat} and~\eqref{eq:delz_bar} we have 
\begin{align}
&\Delta \widehat{Z}_k = B^T(y_k) - Z_k - C + B^T ( \Delta \widehat{y}_k )  && \hspace*{-0.05in} \text{for all} \;\, k \in \mathbb{N}\label{eq:delZ1}\\
&\Delta \overline{Z}_k = B^T( \Delta \overline {y}_k) && \hspace*{-0.4in} \text{for all} \;\, k \in \mathbb{N}. \label{eq:delZ2}
\end{align}
Since all the terms in the right hand side of~\eqref{eq:delZ1} and~\eqref{eq:delZ2} are in $S_{L,M,n}$, then 
\begin{equation}
\Delta \widehat{Z}_k \in S_{L,M,n}, \quad \Delta \overline{Z}_k \in S_{L,M,n}. \label{eq:delZ} 
\end{equation}
Recalling~\eqref{eq:delxhat},~\eqref{Delxbar} and our assumption in~\eqref{eq:assume}, we have 
\begin{equation}
\Delta \overline{X}_k \in S_{L,M,n}. \label{eq:delxbar}
\end{equation}
According to~\eqref{eq:delxhat},~\eqref{eq:delZ} and~\eqref{eq:delxbar}, the total step directions are in $S_{L,M,n}$, 
\begin{align}
&\Delta X_k = \Delta \widehat{X}_k + \Delta \overline{X}_k \in S_{L,M,n}\\
&\Delta Z_k = \Delta \widehat{Z}_k + \Delta \overline{Z}_k \in S_{L,M,n},
\end{align}
and it follows that 
\begin{align}
& X_{k+1} = X_k + t_p \Delta X_k \in S_{L,M,n}\\
& Z_{k+1} = Z_k + t_p \Delta Z_k \in S_{L,M,n}.
\end{align}
Thus, for any $y \in \mathbb{R}^{K}$ and $k \in \mathbb{N} $, if $X_k,Z_k \in S_{L,M,n}$, we have $X_{k+1},Z_{k+1} \in S_{L,M,n}$. Since we have assumed that the initial values $X_0, Z_0 \in S_{L,M,n}$, we conclude by induction that $X_k \in S_{L,M,n}$ and $Z_k \in S_{L,M,n}$ for all $k \in \mathbb{N}$.
\end{proof}

\subsection {Parallel Implementation}
\label{sec:Implementation}

In this section, a parallel algorithm for solving the SDP problems associated with Polya's algorithm is provided. We show how to map the block-diagonal structure of the primal variable and Newton updates described in Section~\ref{sec:Interior-point} to a parallel computing structure consisting of a central root processor with $N$ slave processors. Note that processor steps are simultaneous and transitions between root and processor steps are synchronous. Processors are idle when root is active and vice-versa. A C++ implementation of this algorithm, using MPI and Blas/Lapack libraries is provided at: www.sites.google.com/a/asu.edu/kamyar/software.
Let $N$ be the number of available processors and $J= \mathtt {floor}\left( \frac{L+M}{N} \right)$. As per Algorithm 1, we assume processor $i$ has access to the sub-blocks $ \mathbf{\overline{C}}_i $ and $ \mathbf{\overline{\textbf{B}}}_{j,i} $ defined in~\eqref{eq:Cbar} and~\eqref{eq:Bbar} for $j=1, \cdots, K$. Be aware that minor parts of Algorithm~\ref{alg:SDP} have been abridged in order to simplify the presentation. 

\subsection {Computational Complexity Analysis: SDP Algorithm}
\label{sec:complexity}

NC $\subset$ P is defined to be the class of problems which can be solved in a poly-logarithmic number of steps using a polynomially number processors and is often considered to be the class of problems that can be parallelized efficiently. The class P-complete is a set of problems which are equivalent up to an NC reduction, but contains no problem in NC and is thought to be the simplest class of "inherently sequential" problems. It has been proven that Linear Programming (LP) is P-complete~\cite{limits} and SDP is P-hard (at least as hard as any P-complete problem) and thus is unlikely to admit a general-purpose parallel solution. Given this fact and given the observation that the problem we are trying to solve is NP-hard, it is important to thoroughly understand the complexity of the algorithms we are proposing and how this complexity scales with various parameters which define the size of the problem. To better understand these issues, we have broken our complexity analysis down into several cases which should be of interest to the control community. Note that the cases below do not discuss memory complexity. This is because in the cases when a sufficient number of processors are available, for a system with $n$ states, the memory requirements per block are simply proportional to $n^2$. \vspace*{0.05in}

 \textit {1) Case 1:} \textit{Systems with large number of states}

Suppose we are considering a problem with $n$ states. For

\begin{algorithm}
\vspace*{0.15in}
\textbf{\textit{Inputs}:} $\overline{\textbf{C}}_i,\overline{\textbf{B}}_{j,i}$ for $i=1,\cdots,N$ and $j=1,\cdots,K$ - the sub-blocks of the SDP elements  provided to processor $i$ by the set-up algorithm. \vspace*{0.1in}

\textbf{\textit{Processors Initialization step:}}\\
\For{$i = 1, \cdots, N$, processor $i$}{
Initialize primal and dual variables $\mathbf{X}^0_i$, $\mathbf{Z}^0_i$ and $y^0$ as 
\begin{equation*}
\mathbf{X}^0_i =
\begin{cases}
    I_{(J +1) n },  &   0 \le i < L+M-N J  \\
	I_{J n },  &   L+M-N J \le  i  < N,
\end{cases},
\end{equation*}
\begin{equation*}
  \mathbf{Z}^0_i = \mathbf{X}^0_i \quad \text{and} \quad y^0 = \vec{0} \in \mathbb{R}^K, \label{eq:init}
\end{equation*}
Calculate the complementary slackness~\cite{helmberg2000spectral} $S_{i} = tr( \mathbf{Z}_i^0 \mathbf{X}_i^0 )$. Send $S_{i}$ to processor root.
} \vspace*{0.15in}

\textbf{\textit{Root Initialization step:}}\\
Root processor \textbf{do}\\ 

\hspace*{0.05in} Calculate the barrier parameter~\cite{helmberg2000spectral} $\mu = \frac{1}{3} \sum_{i=1}^{N} S_{i}$.
\hspace*{0.05in} Set the SDP element $ a= \vec{1} \in \mathbb{R}^K $.\\ \vspace*{0.15in}

\textbf{\textit{Processors step 1:}}\\ 
\For{$i=1, \cdots,N$, processor $i$}{ 
\For{$k=1, \cdots,K$}{
Calculate the elements of $\Omega_1$ (R-H-S of~\eqref{eq:SAE1}) 
\begin{small}
\[
\omega_{i,k}  =  tr \left( \mathbf{\overline{\textbf{B}}}_{k,i} (\mathbf{{Z}}_i)^{-1} \left(-\sum_{j=1}^K y_j \mathbf{\overline{\textbf{B}}}_{j,i} + \mathbf {Z}_i + \mathbf{\overline{C}}_i \right) \mathbf{X}_i \right) \label{eq:T1}  
\]
\end{small}
\For{$l=1, \cdots,K$}{
Calculate the elements of the SCM as
\begin{equation}
 \hspace*{-0.3in} \lambda_{i,k,l} =  tr \left( \mathbf{\overline{\textbf{B}}}_{k,i} (\mathbf{{Z}}_i)^{-1} \mathbf{\overline{\textbf{B}}}_{l,i} \mathbf{X}_i \right) \label{eq:T2} 
\end{equation}
}}
Send $\omega_{i,k}$ and $\lambda_{i,k,l}$, $k=1, \cdots,K$ and $l=1, \cdots,K$ to root processor.
} \vspace*{0.15in}

\textbf{\textit{Root step 1:}}\\ 
Root processor \textbf{do} \\ 
\hspace*{0.05in} Construct the R-H-S of~\eqref{eq:SAE1} and the SCM as 
\begin{equation*}
\hspace*{0.1in} \Omega_1 = \left( \begin{array}{ccc}
  \sum_{i=1}^{N} \omega_{i,1} \\
  \sum_{i=1}^{N} \omega_{i,2} \\
  \vdots\\
  \sum_{i=1}^{N} \omega_{i,K}
\end{array} \right)-a \quad \text{and}
\end{equation*}
\begin{equation*}
\Lambda= \left[ \left(\begin{array}{ccc}
  \sum_{i=1}^{N} \lambda_{i,1,1} \\
  \sum_{i=1}^{N} \lambda_{i,2,1} \\
  \vdots\\
  \sum_{i=1}^{N} \lambda_{i,K,1} \end{array} \right) , \cdots,
  \left( \begin{array}{ccc}
  \sum_{i=1}^{N} \lambda_{i,1,K} \\
  \sum_{i=1}^{N} \lambda_{i,2,K} \\
  \vdots \\
  \sum_{i=1}^{N} \lambda_{i,K,K} \end{array}
  \right) \right]
\end{equation*}
\hspace*{0.1in} Solve the following system of equations for the predictor dual step ${\Delta \widehat{y}} \in \mathbb{R}^K $ and send \hspace*{0.1in} ${\Delta \widehat{y}}$ to all processors. 
\begin{equation}
 \hspace*{0.05in} \Lambda {\Delta \widehat{y}} = \Omega_1   \label{eq:SAE1}
\end{equation}

\caption{The parallel SDP solver algorithm}
\label{alg:SDP}
\end{algorithm}

\begin{algorithm}
 \textbf{\textit{Processors step 2:}}\\ 
\For{$i=1, \cdots,N$, processor $i$}{
Calculate the predictor step directions 
\begin{align*}
 \Delta \mathbf{\widehat{X}}_i& = -\mathbf{X}_i \\ 
 + (&\mathbf{{Z}}_i)^{-1} \hspace*{-0.05in} \left( -\sum_{j=1}^K y_j \mathbf{\overline{\textbf{B}}}_{j,i}  + \mathbf {Z}_i + \mathbf{\overline{C}}_i \right) \hspace*{-0.05in} \sum_{j=1}^K \Delta \widehat{y}_j \, \mathbf{\overline{\textbf{B}}}_{j,i} \; \mathbf{X}_i,
\end{align*}
\[
\Delta \mathbf{\widehat{Z}}_i = \sum_{j=1}^K y_j \mathbf{\overline{\textbf{B}}}_{j,i} - \mathbf{Z}_i - \mathbf{\overline{C}}_i + \sum_{j=1}^K \Delta \widehat{y}_j \mathbf{\overline{\textbf{B}}}_{j,i}.
\]
\For{$k=1, \cdots,K$}{
Calculate the elements of $\Omega_2$ (R-H-S of~\eqref{eq:SAE2}) 
\[
 \hspace*{-0.1in}\delta_{i,k} = tr(\mathbf{\overline{\textbf{B}}}_{k,i} (\mathbf{{Z}}_i)^{-1})  ,
 \tau_{i,k} = tr(\mathbf{\overline{\textbf{B}}}_{k,i} (\mathbf{{Z}}_i)^{-1} \Delta \mathbf{\widehat{Z}}_i \Delta \mathbf{ \widehat{X}}_i) 
\]
}
Send $\delta_{i,k}$ and $\tau_{i,k}$, $k=1, \cdots,K$ to root processor.
} \vspace*{0.1in}

 \textbf{\textbf{\textit{Root step 2:}}}\\ 
Root processor \textbf{do}\\
\hspace*{0.05in} Construct the R-H-S of~\eqref{eq:SAE2} as \vspace*{-0.05in}
\begin{align*}
 \Omega_2 =
 \mu
 &\begin{bmatrix}
\sum_{i=1}^{N} \delta_{i,1} & \sum_{i=1}^{N} \delta_{i,2} & \cdots & \sum_{i=1}^{N} \delta_{i,K}
\end{bmatrix}^T
- \\
  &\begin{bmatrix}
  \sum_{i=1}^{N} \tau_{i,1} & \sum_{i=1}^{N} \tau_{i,2} &
  \cdots &
  \sum_{i=1}^{N} \tau_{i,K} \end{bmatrix}^T  
\end{align*}

\hspace*{0.05in} Solve the following system of equations for the corrector dual variable $\Delta \overline{y}$ and send $\Delta \overline{y}$ \hspace*{0.07in} to all processors. \vspace*{-0.1in}
\begin{equation}
\Lambda \Delta \overline{y} = \Omega_2  \label{eq:SAE2}
\end{equation}

\textbf{\textit{Processors step 3:}}\\ 
\For{$i=1, \cdots,N$, processor $i$}{
Calculate the corrector step directions as follows. \vspace*{-0.05in}
\begin{equation*}
\Delta \mathbf{ \overline{Z}}_i = \sum_{j=1}^K \Delta \overline {y}_j \mathbf{\overline{\textbf{B}}}_{j,i}  \vspace*{-0.05in}
\end{equation*}
\begin{equation*}
\Delta \mathbf{ \overline{X}}_i = -(\mathbf{{Z}}_i)^{-1} ( \Delta \mathbf{ \overline{Z}}_i \mathbf{X}_i + \Delta \mathbf{ \widehat{Z}}_i \Delta \mathbf{ \widehat{X}}_i) + \mu (\mathbf{{Z}}_i)^{-1}  \vspace*{-0.05in}
\end{equation*}

Calculate primal dual step total step directions as follows.  \vspace*{-0.05in}
\begin{equation*}
\hspace*{-0.1in}\Delta \mathbf{ X}_i  = \Delta \mathbf{\widehat{X}}_i + \Delta \mathbf{ \overline{X}}_i,
\Delta \mathbf{ Z}_i = \Delta \mathbf{\widehat{Z}}_i + \Delta \mathbf{ \overline{Z}}_i,
\Delta {y} = \Delta \widehat {y} + \Delta \overline {y}.   \vspace*{-0.05in}
\end{equation*}
Set primal step size $t_p$ and  dual step size $t_d$ using an appropriate line search methos.\\ 
 Update primal and dual variables as 
\[
\mathbf{X}_i \equiv \mathbf{X}_i + t_p \Delta \mathbf{X}_i, \quad \mathbf{Z}_i \equiv \mathbf{Z}_i + t_d \Delta \mathbf{Z}_i , \quad y \equiv y  + t_d \Delta {y} \vspace*{-0.15in}
\]
} \vspace*{0.1in}

\textbf{\textit{Processors step 4:}}\\ 
\For{$i=1, \cdots,N$, processor $i$}{
Calculate the contribution to primal cost $ \tilde{\phi}_{i}  = tr \left( \mathbf{\overline{C}}_i \mathbf{X}_i \right)$ and the complementary slack $S_{i}  = tr \left( \mathbf{Z}_i \mathbf{X}_i \right) $.
Send $S_{i}$ and $\tilde{\phi}_{i}$ to root processor.
} \vspace*{0.1in}

\textbf{\textit{Root step 4:}}\\ 
Root processor \textbf{do} \\ 
\hspace*{0.05in} Update the barrier parameter $\mu = \frac{1}{3} \sum_{i=1}^{N} S_{i}$. Calculate primal and dual costs as \hspace*{0.05in} $\phi = \sum_{i=1}^{N} \tilde{\phi}_{i} \; \text{and} \; \psi= a^Ty $. If $ | \phi - \psi | > \eps$, then go to Processors step 1; Otherwise  \hspace*{0.05in} calculate the coefficients of \underline{$P(\alpha)$ as $P_i= \sum_{j=1}^{\tilde{N}} E_j y_{(j + \tilde{N}i-1))}$ for $i=1, \cdots, L_0$.}
\end{algorithm}

\noindent this case, the most expensive part of the algorithm is the calculation of the Schur complement matrix $\Lambda$ by the processors in Processors step 1 (and summed by the root in Root step 1, although we neglect this part). In particular, the computational complexity of the algorithm is determined by the number of operations required to calculate~\eqref{eq:T2}, restated here. \vspace*{-0.05in}
\begin{equation}
\lambda_{i,k,l} =  tr \left( \mathbf{\overline{\textbf{B}}}_{k,i} (\mathbf{{Z}}_i)^{-1} \mathbf{\overline{\textbf{B}}}_{l,i} \mathbf{X}_i \right) \vspace*{-0.05in}
\end{equation}
\begin{equation*} 
 \text{for} \quad k= 1, \cdots, K \; \text{and} \quad l=1, \cdots, K. \vspace*{-0.05in}
\end{equation*}
Since the cost of $n \times n$ matrix-matrix multiplication is $\sim n^3$ and each of $\mathbf{X}_i, \mathbf{{Z}}_i, \mathbf{\overline{\textbf{B}}}_{l,i}$ has $\mathtt{floor}(\frac{L+M}{N})$ number of blocks in $\mathbb{R}^{n \times n}$, the number of operations performed by the $i^{th}$ processor to calculate  $\lambda_{i,k,l}$ for $k= 1, \cdots, K$ and $l=1, \cdots, K$ is \vspace*{-0.05in}
\begin{equation}
\begin{cases}
\sim \mathtt{floor} \left( \dfrac{L+M}{N} \right)  K^2 n^3 \; & N < L+M \\
\sim K^2 n^3 \; & N \geq L+M
\end{cases} \label{eq:flop1}
\end{equation}
at each iteration, where $i=1, \cdots,N$. By substituting $K$ in~\eqref{eq:flop1} from~\eqref{eq:K}, for $N \geq L+M$, each processor performs 
\begin{equation}
\sim \dfrac{((d_p+l-1)!)^2}{(d_p!)^2((l-1)!)^2}n^7 \label{eq:flop_proc} 
\end{equation}
operations per iteration. Therefore, for systems with large $n$ and fixed $d_p$ and $l$, the number of operations per processor required to solve the SDP associated with parameter-dependent feasibility problem
$
A(\alpha)^TP(\alpha)+P(\alpha)A(\alpha) \prec 0,
$
is proportional to $n^7$. Solving the LMI associated with the parameter-independent problem
$
\mathrm{A}^TP+P\mathrm{A} \prec 0
$
using our algorithm or most of the SDP solvers such as~\cite{sedumi,csdp,sdpara} also requires $ O(n^7)$ operations per processor. Therefore, if we have a sufficient number of processors, the proposed algorithm solves both the stability and robust stability problems by performing $O(n^7)$ operations per processor in this case.

\vspace*{0.05in}
 \textit {2) Case 2:} \textit{High Accuracy/Low Conservativity}

In this case we consider the effect of raising Polya's exponent. Consider the definition of simplex as follows. \vspace*{-0.05in}
\begin{equation}
\tilde{\Delta}_l=\left\lbrace \alpha\in \mathbb{R}^l , \sum_{i=1}^{l} \alpha_i=r, \alpha_i\geqslant 0 \right\rbrace \vspace*{-0.05in}
\end{equation}
Suppose we now define the accuracy of the algorithm as the largest value of $r$ found by the algorithm (if it exists) such that if the uncertain parameters lie inside the corresponding simplex, the stability of the system is verified. Typically, increasing Polya's exponent $d$ in~\eqref{eq:polya's exponent} improves the accuracy of the algorithm. If we again only consider Processor step 1, according to~\eqref{eq:flop_proc}, the number of processor operations is independent of the Polya's exponent $d_1$ and $d_2$! Because this part of the algorithm does not vary with Polya's exponent, we look at the root processing requirements associated with solving the systems of equations in~\eqref{eq:SAE1} and~\eqref{eq:SAE2} in Root step 1 using Cholesky factorization. Each of these systems consists of $K$ equations. The computational complexity of Cholesky factorization is $O(K^3)$. Thus the number of operations performed by the root processor is proportional to \vspace*{-0.05in}
\begin{equation} 
K^3 = \dfrac{((d_p+l-1)!)^3}{(d_p!)^3((l-1)!)^3}n^6. \label{eq:flop_root}\vspace*{-0.05in}
\end{equation}
In terms of communication complexity, the most significant operation between the root and other processors is sending and receiving  $\lambda_{i,k,l}$ for $i=1, \cdots, N$, $k=1, \cdots,K$ and $l=1, \cdots,K$ in Processors step 1 and Root step 1. Thus the total communication cost for $N$ processors per iteration is \vspace*{-0.05in}
\begin{equation}
\sim N \cdot K^2 = N \dfrac{((d_p+l-1)!)^2}{(d_p!)^2((l-1)!)^2}n^4. \vspace*{-0.05in} \label{eq:com_cost}
\end{equation}
From~\eqref{eq:flop_proc},~\eqref{eq:flop_root} and~\eqref{eq:com_cost} it is observed that the number of processors operations, root operations and communication operations are independent of Polya's exponent $d_1$ and $d_2$. Therefore, we conclude that for a fixed $d_p$ and sufficiently large number of processors $N$ ($N \geq  L+M$), improving the accuracy by increasing $d_1$ and $d_2$ does not add any computation per processor or communication overhead. \vspace*{0.1in}

 \textit {3) Case 3:} \textit{Algorithm scalability/Speed-up}

\begin{figure}[t]
   \centering
   \includegraphics[scale=0.35]{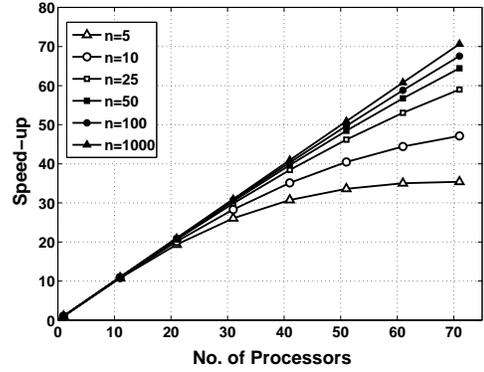}
   \vspace*{-0.1in}
   \caption{Theoretical speed-up vs. No. of processors for different system dimensions $n$ for $l=10$, $d_p=2$, $d_a=3$
and $d_1=d_2=4$, where $L+M=53625$}
   \label{fig:theo_speedup} 
\end{figure}

The speed-up of a parallel algorithm is defined as $\textit{\text{SP}}_N=\dfrac{T_s}{T_N}, $
where $T_s$ is the execution time of the sequential algorithm and $T_N$ is the execution time of the parallel algorithm using $N$ processors. The speed-up is governed by
\vspace*{-0.05in}
\begin{equation}
\textit{\text{SP}}_N=\dfrac{N}{D+NS}, \label{eq:speedup} \vspace*{-0.05in}
\end{equation}
where $D$ is defined as the ratio of the total operations performed by all processors except root to total operations performed by all processors and root. $S$ is the ratio of the operations performed by root to total operations performed by all processors and root. Suppose that the number of available processors is equal to the number of sub-blocks in $C$ defined in~\eqref{eq:C}. Using the above definitions for $D$ and $S$, Equation~\eqref{eq:flop_proc} as the decentralized computation and~\eqref{eq:flop_root} as the centralized computation, $D$ and $S$ can be approximated as
\begin{equation}
D \simeq \dfrac{N \dfrac{((d_p+l-1)!)^2}{(d_p!)^2((l-1)!)^2}n^7}{N \dfrac{((d_p+l-1)!)^2}{(d_p!)^2((l-1)!)^2}n^7+\dfrac{((d_p+l-1)!)^3}{(d_p!)^3((l-1)!)^3}n^6} \;\;
\text{and}
\end{equation}
\begin{equation}
S \simeq \dfrac{\dfrac{((d_p+l-1)!)^3}{(d_p!)^3((l-1)!)^3}n^6}{N \dfrac{((d_p+l-1)!)^2}{(d_p!)^2((l-1)!)^2}n^7+\dfrac{((d_p+l-1)!)^3}{(d_p!)^3((l-1)!)^3}n^6}.
\end{equation}
According to~\eqref{eq:L} and~\eqref{eq:M} the number of processors $N=L+M$ is independent of $n$; Therefore
\vspace*{-0.05in}
\[
\lim_{n\to\infty} D = 1 \quad \text{and} \quad \lim_{n\to\infty} S = 0. \vspace*{-0.05in}
\]
By substituting $D$ and $S$ in~\eqref{eq:speedup} with their limit values, we have $\lim_{n\to\infty} \textit{\text{SP}}_N = N$. Thus, for large $n$, by using $L+M$ processors the presented decentralized algorithm solves large robust stability problems $L+M$ times faster than the sequential algorithms. For different values of the state-space dimension $n$, the theoretical speed-up of the algorithm versus the number of processors is illustrated in Fig.~\ref{fig:theo_speedup}. As shown in Fig.~\ref{fig:theo_speedup}, for problems with large $n$, by using $N \leq L+M$ processors the parallel algorithm solves the robust stability problems approximately $N$ times faster than the sequential algorithm. As $n$ increases, the trend of speed-up becomes increasingly linear. Therefore, in case of problems with a large number of states $n$, our algorithm becomes increasingly efficient in terms of processor utilization.

 \textit {4) Case 4:} \textit{Synchronization and load balancing}

\begin{figure}[t]
   \centering
 \hspace*{-0.65in}  \includegraphics[scale=0.3]{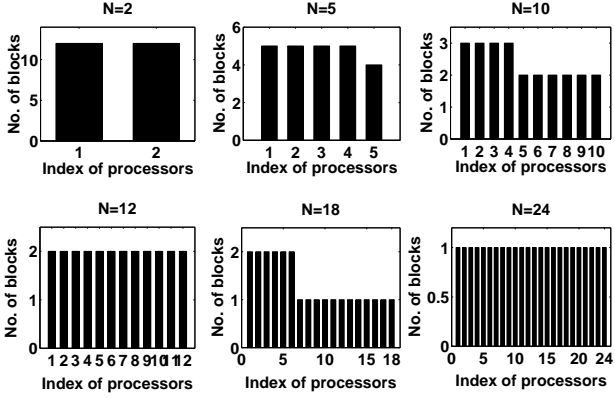} \vspace*{-0.1in}
   \caption{The number of blocks of the SDP elements assigned to each processor; An illustration of load balancing} 
   \label{fig:load_balance} 
\end{figure}

The proposed algorithm is synchronous in that all processors must return values before the centralized step can proceed. However, in the case where we have fewer processors than blocks, some processors may be assigned one block more than other processors. In this case, some processors may remain idle while waiting for the more heavily loaded blocks to complete. In the worst case. this can result in a 50\% decrease in speed. We have addressed this issue in the following manner:
\begin{enumerate}
\item We allocate almost the same number ($\pm 1$) of blocks of the SDP elements $C$ and $B_{i}$ to all processors, i.e., $\mathtt{floor}(\frac{L+M}{N})+1$ blocks to $r$ processors and $\mathtt{floor}(\frac{L+M}{N})$ blocks to the other $N-r$ processors, where $r$ is the remainder of dividing $L+M$ by $N$.
\item We assign the same routine to all of the processors in the Processors steps of Alg. 2.
\end{enumerate}
If $L+M$ is a multiple of $N$, then the algorithm assigns the same amount of data, i.e., $\frac{L+M}{N}$ blocks of $C$ and $B_{i}$ to each processor. In this case, the processors are perfectly synchronized. If $L+M$ is not a multiple of $N$, then according to~\eqref{eq:flop1}, $r$ of $N$ processors perform $K^2n^3$ extra operations per iteration. This fraction is $\dfrac{1}{1+\mathtt{floor}(\frac{L+M}{N})} \leq 0.5$ of the operations per iteration performed by each of $r$ processors. Thus in the worst case, we have a 50\% reduction, although this situation is rare. As an example, the load balancing (distribution of data and calculation) for the case of solving an SDP of the size $L+M=24$ using different numbers of available processors $N$ is demonstrated in Fig.~\ref{fig:load_balance}. This figure shows the number of blocks that are allocated to each processor. According to this figure, for $N=2,12$ and 24, the processors are well-balanced, whereas for the case where $N=18$, twelve processors perform 50$\%$ fewer calculations.

  \textit{5) Case 5:} \textit{Communication graph}

\begin{figure}[t]
\centering
\includegraphics[scale=0.5]{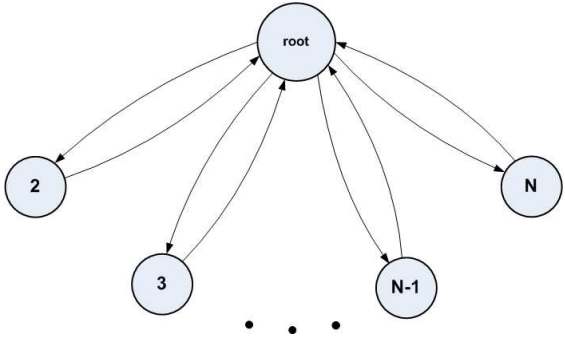} \vspace*{-0.05in}
\caption{The communication graph of the SDP algorithm}
\label{fig:SDP_graph}  \vspace*{-0.15in}
\end{figure}

The communication directed graph of the SDP algorithm (Fig.~\ref{fig:SDP_graph}) is static (fixed for all iterations). At each iteration, root sends messages ($\Delta \widehat{y}$ and $\Delta \overline{y}$) to all of the processors and receives messages ($\lambda_{i,k,l}$ in~\eqref{eq:T2}) from all of the processors. The adjacency matrix of the communication directed graph is defined as follows. For $i=1,\cdots,N$ and $j=1, \cdots, N$,
\[
[T_G]_{i,j} :=
\begin{cases}
1 \quad & \text{if} \; \big( i=1 \; \text{or} \; j=1 \big) \; \text{and} \; \big( i \neq j \big)     \\
0 \quad & \text{Otherwise}.
\end{cases}   \vspace*{-0.2in}
\]

\section{TESTING AND VALIDATION}
\label{sec:RESULTS} \vspace*{-0.05in}

In this section, we present validation data in 4 key areas. First, we present analysis results for a realistic large-scale model of tokamak operation using a discretized PDE model. Next we present accuracy and convergence data and compare our algorithm to the SOS approach. Next, we analyze scalability and speed-up of our algorithm as we increase the number of processors and compare our results to the general-purpose parallel SDP solver SDPARA. Finally, we explore the limits of the algorithm in terms of problems size when implemented on a moderately powerful cluster computer and using a moderate processor allocation on the Blue Gene supercomputer. \vspace*{0.05in}

\textit{1) Example 1:} \textit{Application to control of a discretized PDE model in fusion research.}

The goal of this example is to use the proposed algorithm to solve a real-world stability problem. A simplified model for the poloidal magnetic flux gradient in a Tokamak reactor~\cite{Tokamak} is  \vspace*{-0.05in}
\begin{equation}
\dfrac{\partial \psi_x (x,t)}{\partial t} = \dfrac{1}{\mu_0 a^2} \dfrac{\partial}{\partial x} \left( \dfrac{\eta(x)}{x} \dfrac{\partial}{\partial x} \left( x \psi_x (x,t) \right) \right) 
\end{equation}
with the boundary conditions $\psi_x(0,t)=0$ and $\psi_x(1,t)=0$, where $\psi_x$ is the deviation of the flux gradient from a reference flux gradient profile, $\mu_0$ is the permeability of free space, $\eta(x)$ is the plasma resistivity and $a$ is the radius of the last closed magnetic surface (LCMS). To obtain the finite-dimensional state-space representation of the PDE, we discretize the PDE in the spatial domain $(0,1)$. The state-space model is then 
\begin{equation}
\dot{\psi}_x(t) = A(\eta(x)) \psi_x(t),  \label{eq:discrete_sys} 
\end{equation}
where $A(\eta(x)) \in \mathbb{R}^{N \times N}$ has the following non-zero entries.
\begin{equation*}
\hspace*{-1.2in} a_{11}= \dfrac{-4}{3  \mu_0 \Delta x^2  a^2} \left( \hspace*{-0.05in} \dfrac{\eta(x_{\frac{3}{2}})}{x_{\frac{3}{2}}} + \dfrac{2 \eta(x_{\frac{3}{4}})}{x_{\frac{3}{4}}} \right),
 \end{equation*}
 \begin{equation}
\hspace*{-1.64in}  a_{12}= \dfrac{4}{3  \mu_0 \Delta x^2  a^2} \left(\dfrac{\eta(x_{\frac{3}{2}})x_2}{x_{\frac{3}{2}}} \right),
\end{equation}
\begin{equation}
\hspace*{-0.04in} a_{j,j-1}= \dfrac{1}{\Delta x^2 \mu_0 a^2} \left( \dfrac{\eta(x_{j-\frac{1}{2}})}{x_{j-\frac{1}{2}}}x_{j-1}  \right)
\; \text{for} \;  j=2, \cdots, N-1
\end{equation}
\begin{align}
 & \hspace*{-0.1in} a_{j,j}=\dfrac{-1}{\Delta x^2 \mu_0 a^2}  \left(\dfrac{\eta(x_{j+\frac{1}{2}})}{x_{j+\frac{1}{2}}}+  \dfrac{\eta(x_{j-\frac{1}{2}})}{x_{j-\frac{1}{2}}}\right)x_j \nonumber \\
& \qquad \qquad \qquad \qquad \qquad \quad \qquad\text{for} \;  j=2, \cdots, N-1 
\end{align}
\begin{equation}
 a_{j,j+1}=  \dfrac{1}{\Delta x^2 \mu_0 a^2} \left( \dfrac{\eta(x_{j+\frac{1}{2}})}{x_{j+\frac{1}{2}}}x_{j+1} \right)  \; \text{for} \;  j=2, \cdots, N-1
\end{equation}
\begin{equation*}
\hspace*{-1.3in} a_{N,N-1}= \dfrac{4}{3 \Delta x \mu_0 a^2} \dfrac{\eta(x_{N-\frac{1}{2}})x_{N-1}}{x_{N-\frac{1}{2}} \Delta x} ,
\end{equation*}
\begin{equation}
\hspace*{-0.1in}  a_{N,N}= \dfrac{-4}{3 \Delta x \mu_0 a^2} \left(\dfrac{2\eta(x_{N+\frac{1}{4}})x_N}{x_{N+\frac{1}{4}} \Delta x} + \dfrac{\eta(x_{N-\frac{1}{2}})x_N}{x_{N-\frac{1}{2}} \Delta x}\right), \vspace*{-0.05in}
\end{equation}
where $\Delta x=\dfrac{1}{N}$ and $x_j:=(j-\frac{1}{2}) \Delta x$.

We discretize the model at $N=7$ points. Typically the $\eta(x_k)$ 
are not precisely known (they depend on other state variables), so we substitute for $\eta(x_k)$ in~\eqref{eq:discrete_sys} with $\widehat{\eta}(x_k)+\alpha_j$, where $\widehat{\eta}(x_k)$ are the nominal values of $\eta(x_k)$ and $\alpha_j$ are the uncertain parameters. At $x_k= 0.036, 0.143, 0.286, 0.429, $ $ 0.571, 0.714, 0.857, 0.964$, we use data from the Tore Supra reactor to estimate the $\widehat{\eta}(x_k)$ as $1.775 \cdot 10^{-8}, 2.703 \cdot 10^{-8}, 5.676 \cdot 10^{-8}, 1.182 \cdot 10^{-7}, 2.058 \cdot 10^{-7}, 3.655 \cdot 10^{-7}, 1.076 \cdot 10^{-6}, 8.419 \cdot 10^{-6}$. The uncertain system is then written as    \vspace*{-0.05in}
\begin{equation}
\dot{\psi}_x(t) = A(\alpha) \psi_x(t), \label{eq:discrete_uncertain} \vspace*{-0.05in}
\end{equation}
where $A$ is affine, $A(\alpha) = A_0 + \sum_{i=1}^8 A_i \alpha_i$ (the $A_i$ are omitted for the sake of brevity).
For a given $\rho$, we restrict the uncertain parameters $\alpha_j$ to $S_\rho$, defined as
\vspace*{-0.07in}
\begin{equation}
S_\rho := \{ \alpha \in \mathbb{R}^8 : \sum_{i=1}^8 \alpha_i = -6|\rho|, -|\rho| \leq \alpha_i \leq |\rho| \}, \vspace*{-0.05in}
\end{equation}
which is a simplex translated to the origin. We would like to determine the maximum value of $\rho$ such that the system is stable by solving the following optimization problem. \vspace*{-0.05in}
\begin{equation*}
\hspace*{-1.85in} \max \quad \rho   \vspace*{-0.05in}
\end{equation*}
\begin{equation}
\hspace*{0.5in} \text{s.t.} \quad \text{System~\eqref{eq:discrete_uncertain} is stable for all} \; \alpha \in S_\rho. \label{eq:optim22} \vspace*{-0.05in}
\end{equation}
To represent $S_\rho$ using the standard unit simplex defined in~\eqref{eq:simplex}, we define the invertible map $g: \Delta_8 \rightarrow S_\rho$ as \vspace*{-0.025in}
\begin{equation}
\hspace*{-0.07in} g(\alpha)=\bmat{g_1(\alpha) & \cdots & g_8(\alpha) },\; g_i(\alpha):= 2|\rho|(\alpha_i-0.5). \vspace*{-0.025in}
\end{equation}
Then, if we let $A'(\alpha) = A(g(\alpha))$, since $g$ is one-to-one, \vspace*{-0.05in}
$$
\{ A(\alpha') :  \alpha' \in S_\rho\} \hspace*{-0.03in} = \hspace*{-0.03in} \{ A(g(\alpha)):  \alpha \in \Delta_8\} \hspace*{-0.03in} =\hspace*{-0.03in} \{ A'(\alpha) :  \alpha \in \Delta_8\}. \vspace*{-0.05in}
$$  
Thus stability of $ \dot{\psi}_x(t) = A'(\alpha) \psi_x(t), \text{ for all }\alpha \in \Delta_l$ is equivalent to stability of Equation~\eqref{eq:discrete_uncertain} for all $\alpha \in S_\rho$.

We solve the optimization problem in~\eqref{eq:optim22} using bisection. For each trial value of $\rho$, we use the proposed parallel SDP solver to solve the associated SDP obtained by the parallel set-up algorithm. The SDP problems have 224 constraints with the primal variable $X \in \mathbb{R}^{1092 \times 1092}$. The normalized maximum value of $\rho$  is found to be $0.0019$. In this particular example, the optimal value of $\rho$ does not change with the degrees of $P(\alpha)$ and Polya's exponents~$d_1$ and $d_2$, primarily because the model is affine.

\begin{figure}[t]
\vspace*{-0.1in}
   \centering
   \includegraphics[scale=0.35]{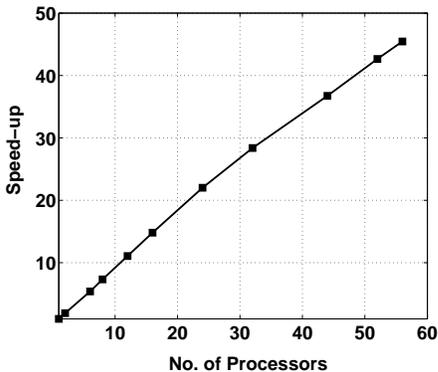}
   \vspace*{-0.05in}
   \caption{Speed-up of set-up and SDP algorithms vs. number of processors for a discretized model of magnetic flux in Tokamak}
      \vspace*{-0.05in}
   \label{fig:Cosmea}   
\end{figure}

The SDPs are constructed and solved on a parallel Linux-based cluster Cosmea at Argonne National Laboratory. Fig.~\ref{fig:Cosmea} shows the algorithm speed-up vs. the number of processors.
Note that solving this problem by SOSTOOLS~\cite{sos1} on the same machine is impossible due to the lack of unallocated memory. 

\textit{2) Example 2:} \textit{Accuracy and Convergence}

\begin{figure}[t]
\centering
  \includegraphics[scale=0.3]{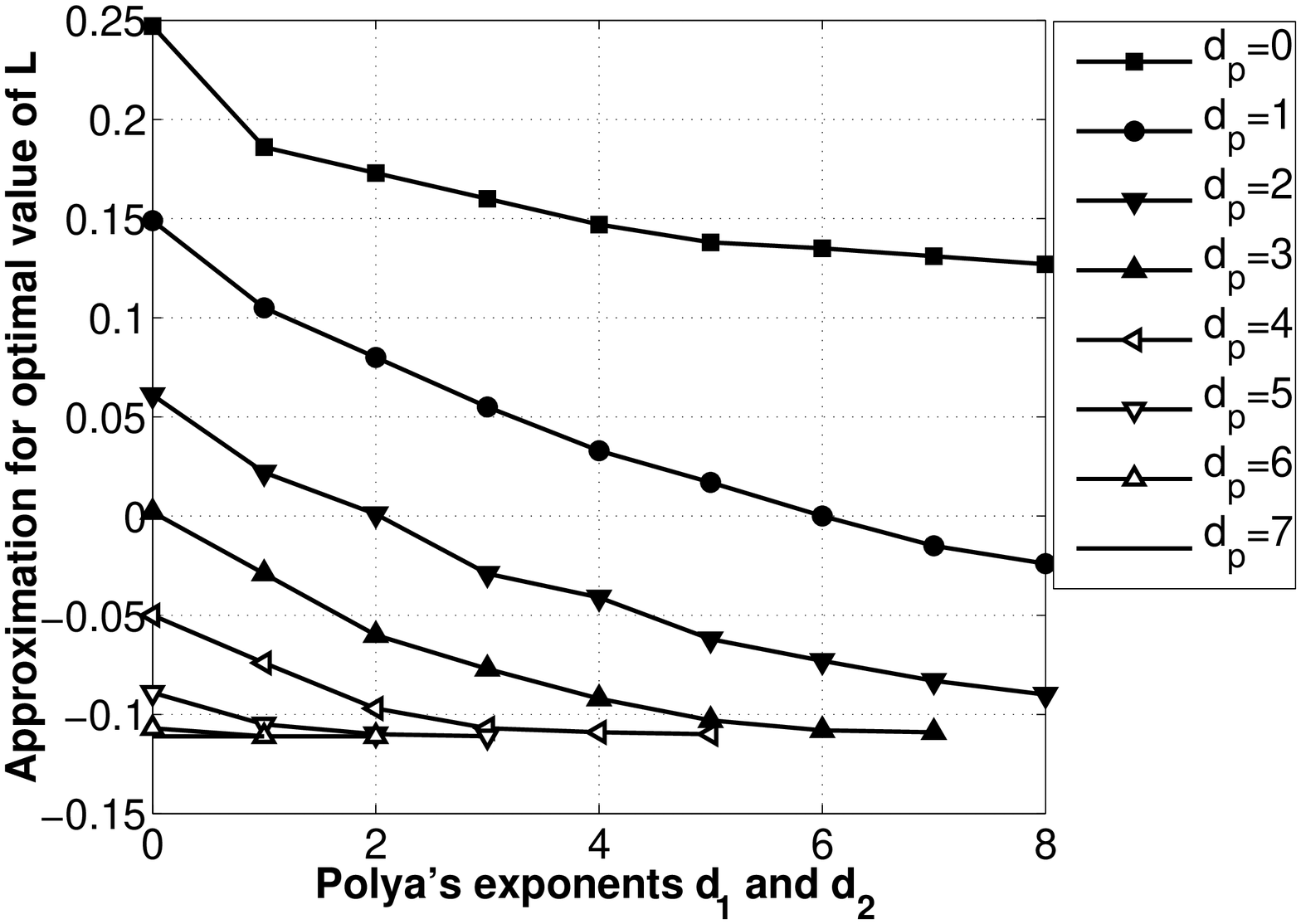} \vspace*{-0.15in}
   \caption{Upper bound on optimal $L$ vs. Polya's exponents $d_1$ and $d_2$, for different degrees of $P(\alpha)$. ($d_1=d_2$).}
   \label{fig:conserve1}
\centering 
  \hspace*{-0.35in} \includegraphics[scale=0.28]{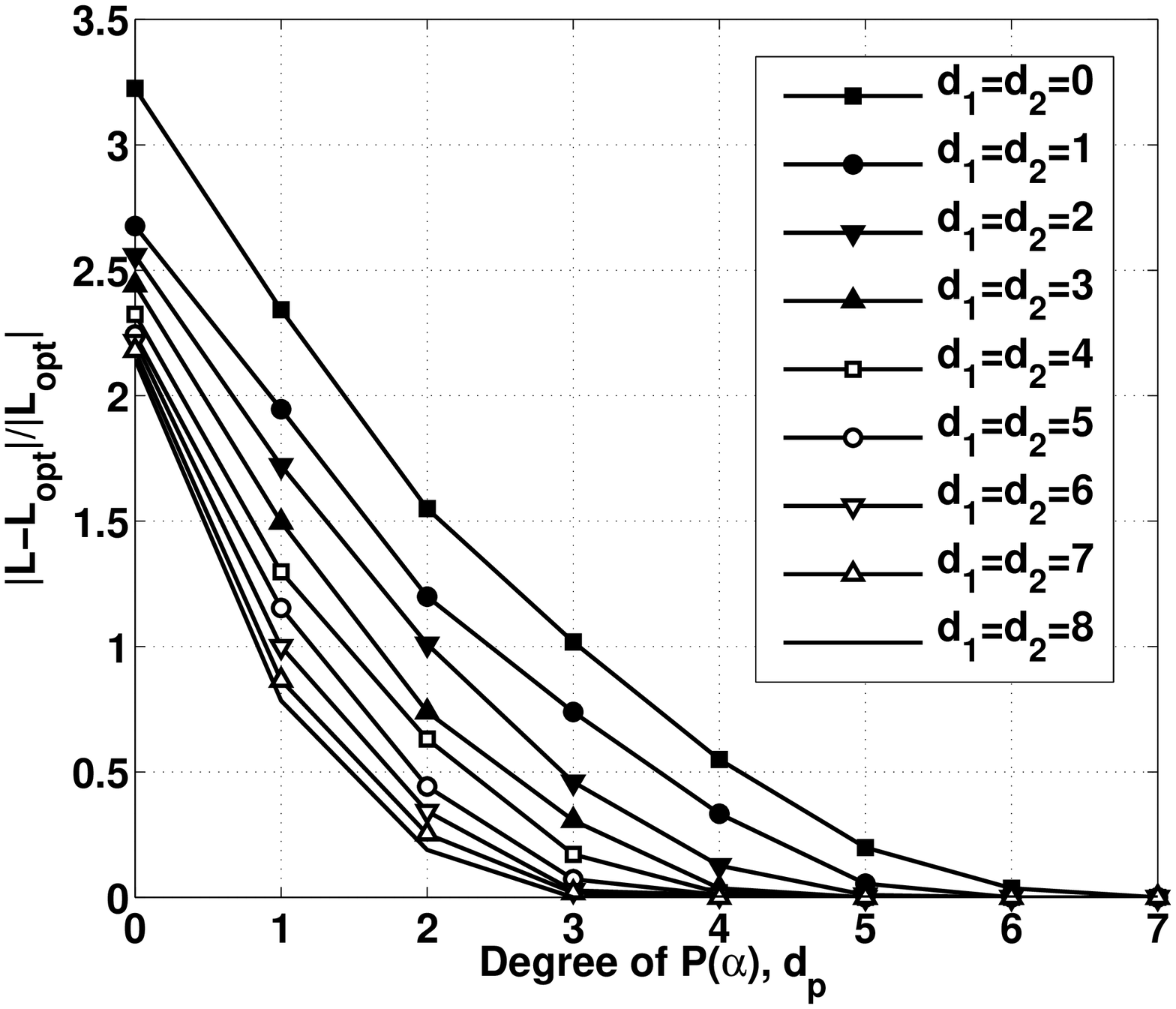} \vspace*{-0.1in}
   \caption{ Error of the approximation for the optimal value of $L$ vs. degrees of $P(\alpha)$, for different Polya's exponents}
   \label{fig:conserve2}   
\end{figure}

The goal of this example is to investigate the effect of the degree of $P(\alpha)$, $d_p$, and the Polya's exponents, $d_1,d_2$ on the accuracy of the algorithm. Given a computer with fixed amount of RAM, we compare the accuracy of the proposed algorithm with SOS algorithm. Consider the system $\dot x(t)=A(\alpha)x(t)$ where $A$ is a polynomial degree 3 defined as \vspace*{-0.05in}
\begin{equation}
A(\alpha)=A_1\alpha_1^3+A_2\alpha_1^2\alpha_2+A_3\alpha_1\alpha_2\alpha_3+A_4\alpha_1\alpha_3^2+A_5\alpha_2^3+A_{6}\alpha_3^3
\vspace*{-0.1in} \end{equation}
with the constraint   \vspace*{-0.05in}
\[
\alpha \in S_L := \left\lbrace \alpha \in \mathbb{R}^3: \sum_{i=1}^3 \alpha_i = 2L+1, L \leq \alpha_i \leq 1 \right\rbrace
\]
\begin{footnotesize}
  \vspace*{-0.1in}
\begin{align*}
&A_1 \hspace*{-0.04in} = \hspace*{-0.04in}
\left[ \hspace*{-0.06in} \begin{array}{ccc}
-0.61 & -0.56 & 0.402 \\
-0.48 & -0.550 & 0.671 \\
-1.01 & -0.918 & 0.029
\end{array} \hspace*{-0.05in} \right],
\hspace*{0.1in} A_2 \hspace*{-0.04in}=\hspace*{-0.06in}
\left[ \hspace*{-0.04in} \begin{array}{ccc}
-0.484 & -0.86 & 1.5 \\
-0.732 & -0.841 & -0.126 \\
 0.685 &  0.305 &  0.106
\end{array} \hspace*{-0.065in} \right] \\
&  A_3 \hspace*{-0.04in}= \hspace*{-0.04in}
\left[ \hspace*{-0.06in} \begin{array}{ccc}
-0.357 &  0.344 & -0.661\\
-0.210 & -0.505 &  0.588\\
0.268 &  0.487 & -0.846 \\
\end{array} \hspace*{-0.06in} \right] \hspace*{-0.03in} , \hspace*{-0.02in}
 A_4 \hspace*{-0.04in}= \hspace*{-0.04in}
\left[ \hspace*{-0.06in} \begin{array}{ccc}
-0.881 & -0.436 & 0.228\\
0.503 & -0.812 & 0.249\\
-0.012 &  0.542 & -0.536
\end{array} \hspace*{-0.06in} \right] \\
& 
 A_5 \hspace*{-0.04in}= \hspace*{-0.04in}
\left[ \hspace*{-0.06in} \begin{array}{ccc}
-0.703 & -0.298 & -0.178\\
0.402 & -0.761 & -0.300\\
-0.010 &  0.461 & -0.588
\end{array} \hspace*{-0.06in} \right] \hspace*{-0.03in} , \hspace*{-0.02in}
 A_6 \hspace*{-0.04in}= \hspace*{-0.04in}
\left[ \hspace*{-0.06in} \begin{array}{ccc}
-0.201 & -0.182 & -0.557\\
0.803 & -0.412 & -0.203\\
-0.440 &  0.011 & -0.881
\end{array} \hspace*{-0.06in} \right]
\end{align*}
\end{footnotesize}
\noindent Defining $g$ as in Example 1, the problem is \vspace*{-0.1in}
\begin{equation*}
\hspace*{-2.6in} \min \quad L \vspace*{-0.1in}
\end{equation*}
\begin{equation}
\hspace*{0.2in} \text{s.t.} \quad \dot{x}(t)=A(g(\alpha))x(t) \; \text{is stable for all} \; \alpha \in \Delta_3. \label{eq:optim} \vspace*{-0.1in}
\end{equation}
Using bisection in $L$, as in Example 1, we varied the parameters $d_p$, $d_1$ and $d_2$. The cluster computer Karlin at the Illinois Institute of Technology with 24 Gbytes/node of RAM (216 Gbytes total memory) was used to run our algorithm. The upper bounds on the optimal $L$ are shown in Fig.~\ref{fig:conserve1} in terms of $d_1$ and $d_2$ and for different $d_p$. Considering the optimal value of $L$ to be $L_{\text{opt}}=-0.111$, Fig.~\ref{fig:conserve1} shows how increasing $d_p$ and/or $d_1,d_2$ - when they are still relatively small - improves the accuracy of the algorithm. Fig.~\ref{fig:conserve2} demonstrates how the error in our upper bound for $L_{\text{opt}}$ decreases by increasing $d_p$ and/or $d_1,d_2$.

For comparison, we solved the same stability problem using the SOS algorithm~\cite{sos1} using only a single node of the same cluster computer and 24 Gbytes of RAM. We used the Positivstellensatz approach based on~\cite{stengle1973nullstellensatz} to impose the constraints $\sum_{i=1}^3 \alpha_i = 2L+1$ and $ L \leq \alpha_i \leq 1$. Table~\ref{tab:SOS} shows the upper bounds on $L$ given by the SOS algorithm using different degrees for $x$ and $\alpha$. By considering a Lyapunov function of degree two in $x$ and degree one in $\alpha$, the SOS algorithm gives $-0.102$ as the upper bound on $L_{opt}$ as compared with our value of $-0.111$. Increasing the degree of $\alpha$ in the Lyapunov function beyond degree two resulted in a failure due to lack of memory. Note that while relevant, this comparison may not be entirely fair as the SOS algorithm has not been decentralized and it can handle global nonlinear stability problems, which our algorithm cannot.

\begin{footnotesize}
\begin{table}[t]
\caption{Upper bounds found for $L_{opt}$ by SOS algorithm using different degrees for $x$ and $\alpha$ (inf: infeasible, O.M.: Out of Memory)} \vspace*{-0.15in}
\label{tab:SOS}
\begin{center}
\begin{tabular}{|c|c|c|c|}
\hline
\backslashbox{  Degree in $x$}{Degree in $\alpha$}
  & 0 & 1 & 2 \\
\hline
1 & inf. & inf. & inf. \\
\hline
2 & inf. & -0.102 & O.M. \\
\hline
3 & inf. & O.M. &  O.M. \\
\hline
\end{tabular} 
\end{center}
\end{table} 
\end{footnotesize}

\textit{3) Example 3:} \textit{Speed-up}

\begin{center}
\begin{figure}[t]
\centering   \hspace*{0.05in}
\hspace*{-0.3in} \includegraphics[scale=0.3]{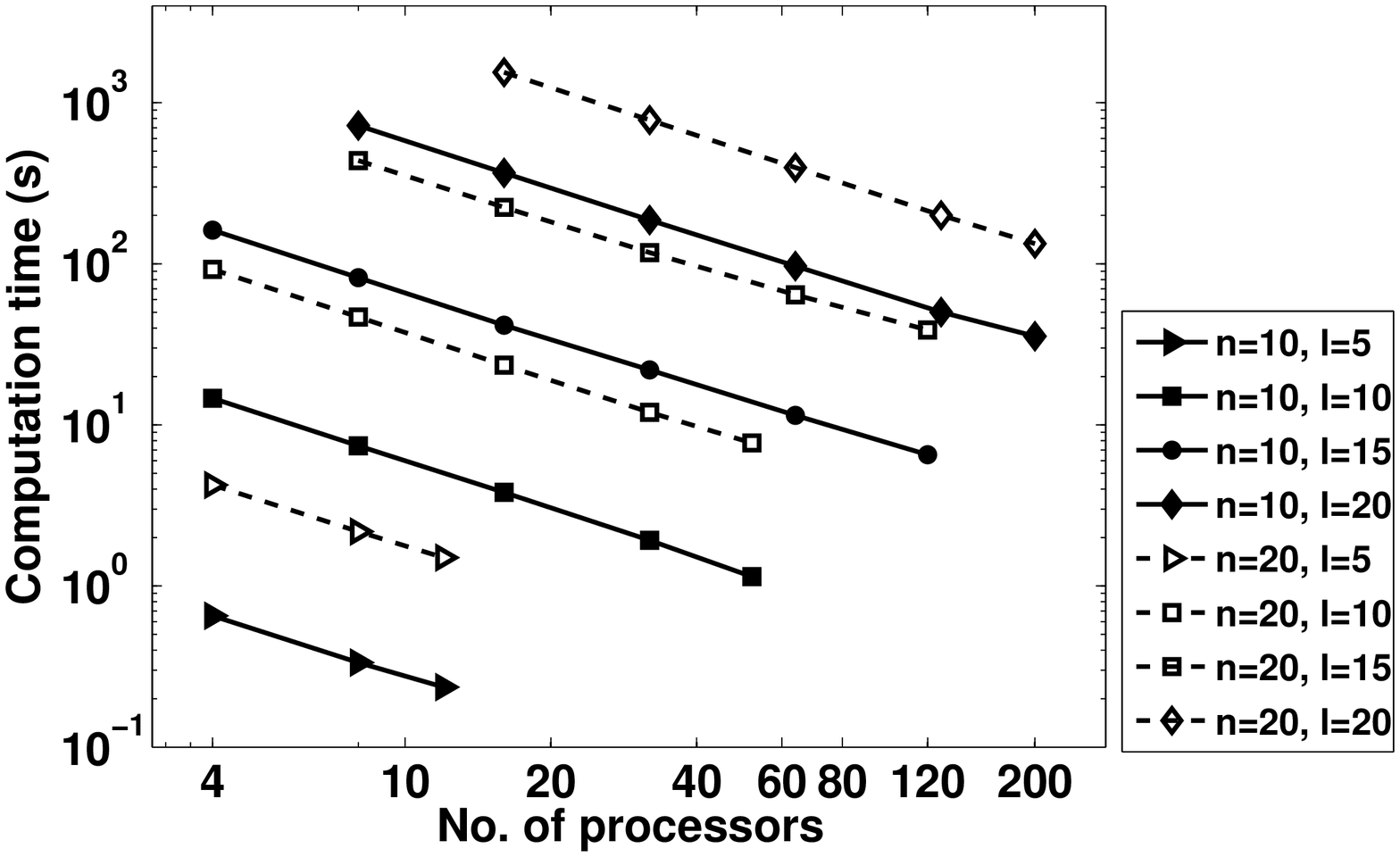} \vspace*{-0.1in} 
 \caption{Computation time of the parallel set-up algorithm vs. number of processors for different dimensions of linear system $n$ and numbers of uncertain parameters $l$- executed on Blue Gene supercomputer of Argonne National Labratory}
\label{fig:speedup1}
\hspace*{-0.35in} \includegraphics[scale=0.3]{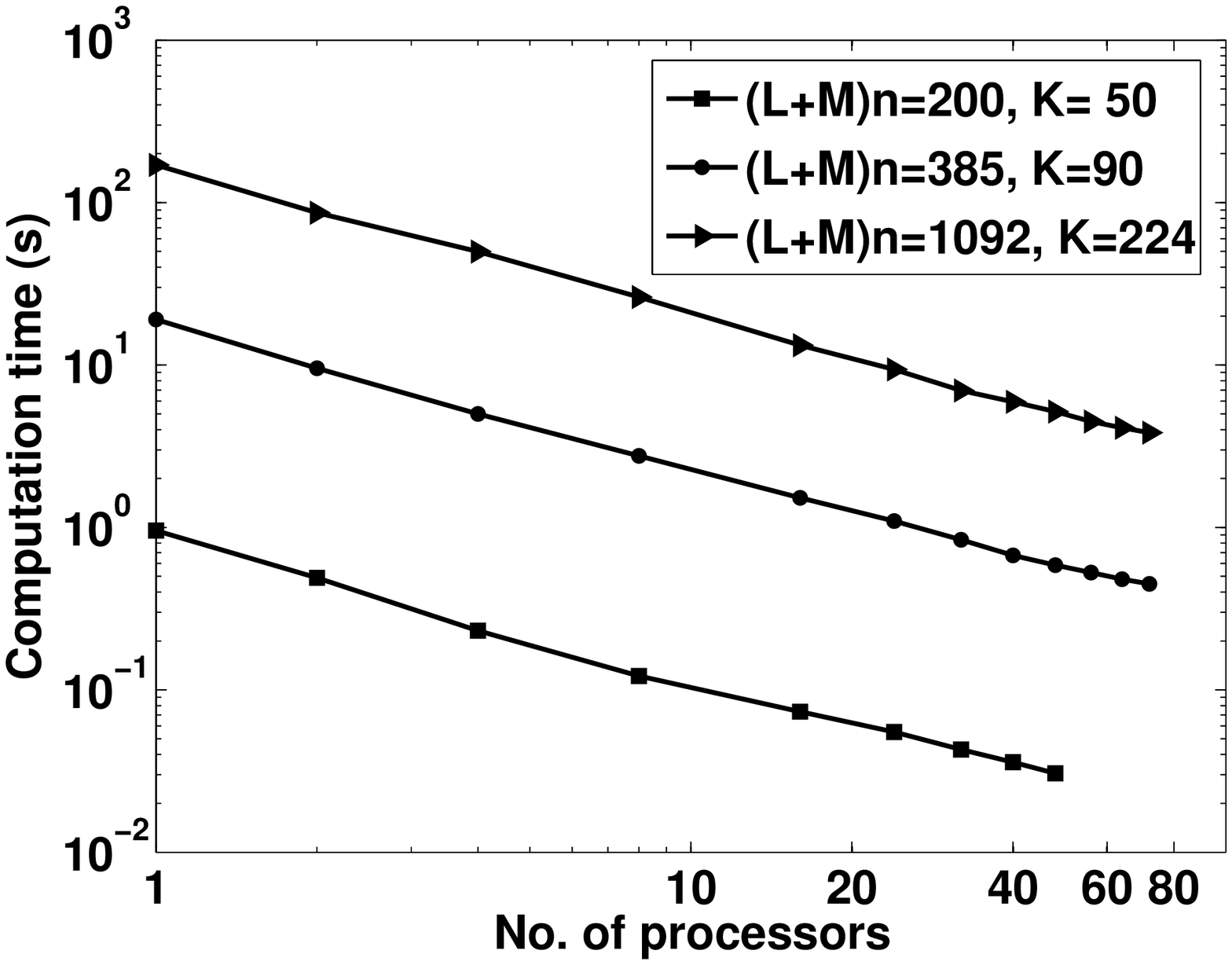} \vspace*{-0.1in}
 \caption{Computation time of the parallel SDP algorithm vs. number of processors for different dimensions of primal variable $(L+M)n$ and of dual variable $K$- executed on Karlin cluster computer of Illinois Institute of Technology}
\label{fig:speedup2} 
\end{figure}
\end{center}

In this example we evaluate the efficiency of the algorithm in using additional processors to decrease computation time. As mentioned in Section~\ref{sec:complexity} on computational complexity, the measure of this efficiency is termed speed-up and in Case 3, we gave a formula for this number. To evaluate the true speed-up, we first ran the set-up algorithm on the Blue Gene supercomputer at Argonne National Laboratory using three random linear systems with different state-space dimensions and numbers of uncertain parameters. Fig.~\ref{fig:speedup1} shows a log-log plot of the computation time of the set-up algorithm vs. the number of processors. As can be seen, the scalability of the algorithm is practically ideal for several different state-space dimensions and numbers of uncertain parameters.

To evaluate the speed-up of the SDP portion of the algorithm, we solved three random SDP problems with different dimensions using the Karlin cluster computer.
Fig.~\ref{fig:speedup2} gives a log-log plot of the computation time of the SDP algorithm vs. the number of processors for three different dimensions of the primal variable $X$ and the dual variable $y$. As indicated in the figure, the three dimensions of the primal variable $X$ are $ 200,\; 385$ and 1092,
 and the dimensions of the dual variable $y$ are $K=50, \; 90$ and 224, respectively.
In all cases, $d_p=2$ and $d_1=d_2=1$. The linearity of the Time vs. Number of Processors curves in all three cases demonstrates the scalability of the SDP algorithm.

\begin{figure}[t]
 \centering
 \includegraphics[scale=0.3]{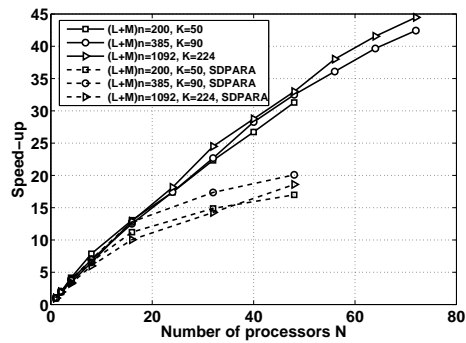} \vspace*{-0.1in}
 \caption{Comparison between the speed-up of the present SDP solver and SDPARA 7.3.1, executed on Karlin cluster computer}
 \label{fig:speedup_sdpara} \vspace*{-0.2in}
 \end{figure}

For comparison, we plot the speed-up of our algorithm vs. that of the general-purpose parallel SDP solver SDPARA 7.3.1 as illustrated in Fig.~\ref{fig:speedup_sdpara}. Although similar for a small number of processors, for a larger number of processors, SDPARA saturates, while our algorithm remains approximately linear.

\vspace*{-0.2in}

\textit{4) Example 4:} \textit{Max state-space and parameter dimensions for a 9-node Linux cluster computer}

\begin{figure}[t] 
\centering
 \includegraphics[scale=0.33]{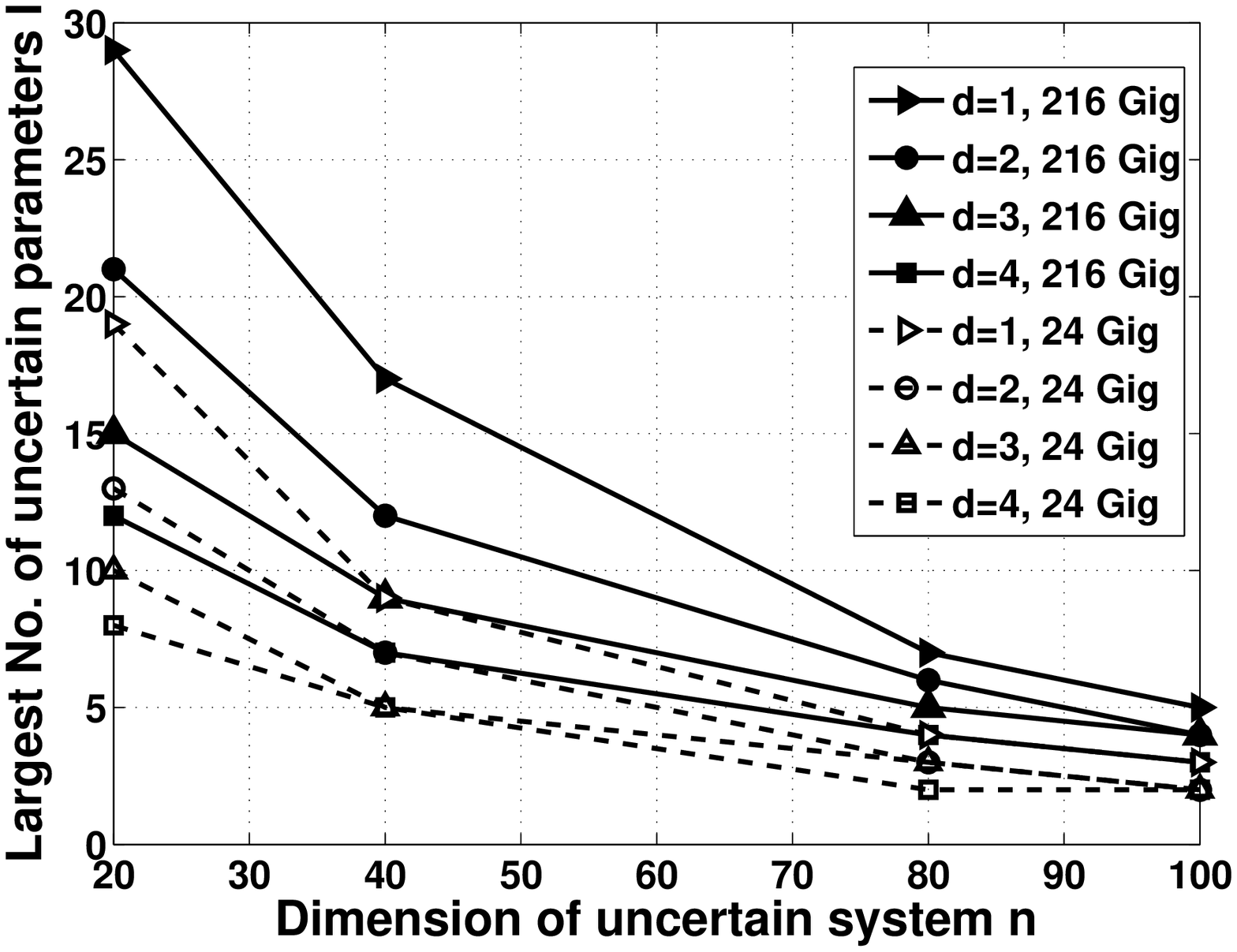}
\centering
\includegraphics[scale=0.33]{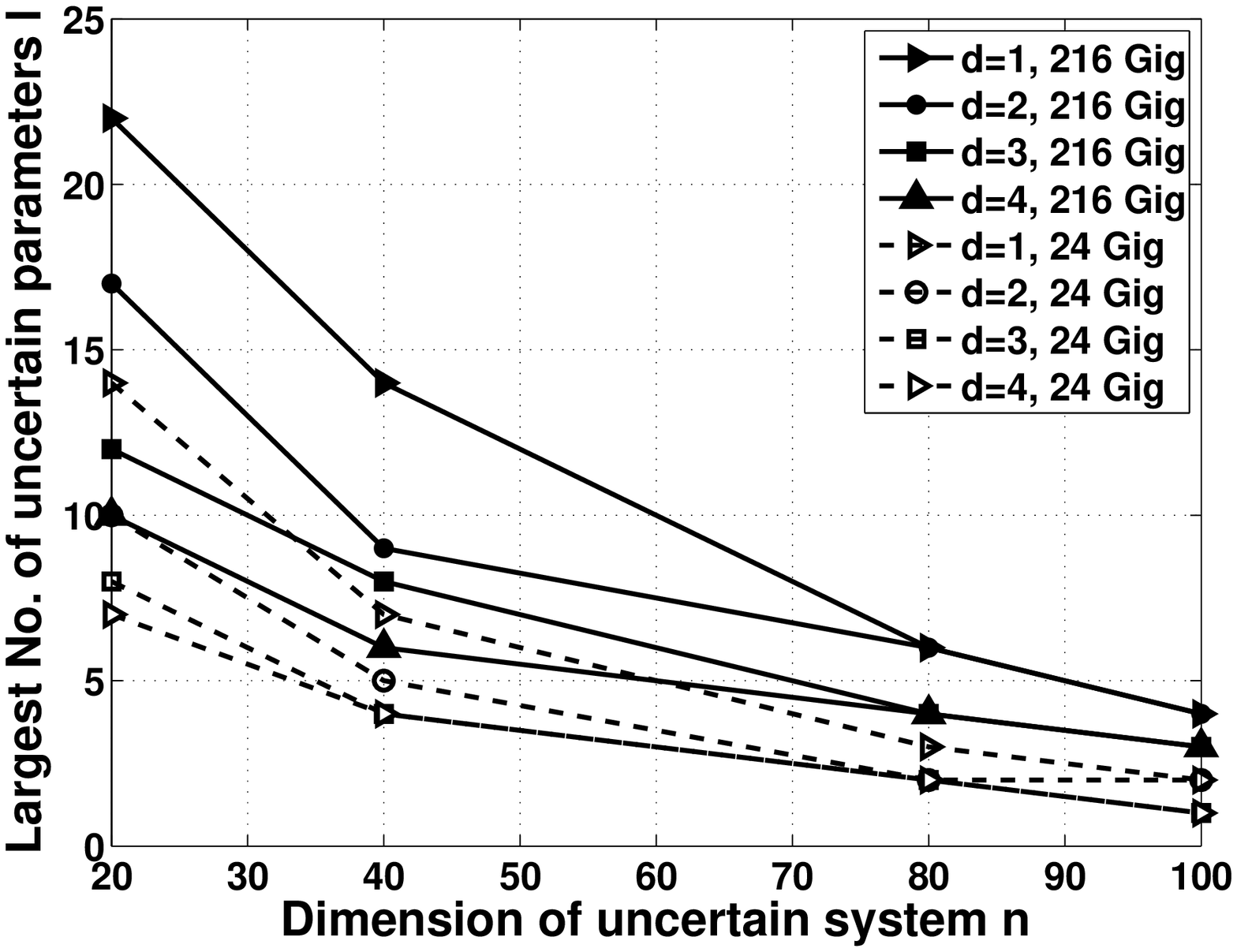}
\caption{Largest number of uncertain parameters of $n$-dimensional systems for which the set-up algorithm (Left) and SDP solver (Right) can solve the robust stability problem of the system using 24 and 216 Gig of RAM}
\label{fig:size_setup}
\end{figure}

The goal of this example is to show that given moderate computational resources, the proposed decentralized algorithms can solve robust stability problems for systems with 100+ states. We used the Karlin cluster computer with 24 Gbytes/node RAM and nine nodes. We ran the set-up and SDP algorithms to solve the robust stability problem with dimension $n$ and $l$ uncertain parameters on one and nine nodes of Karlin cluster computer. Thus the total memory access was thus 24 Gig and 216 Gig, respectively.
Using trial and error, for different $n$ and $d_1,d_2$ we found the largest $l$ for which the algorithms do not terminate due to insufficient memory (Fig.~\ref{fig:size_setup}). In all of the runs $d_a=d_p=1$. Fig.~\ref{fig:size_setup} shows that by using 216 Gbytes of RAM, the algorithms can solve the stability problem of size $n=100$ with 4 uncertain parameters in $d_1=d_2=1$ Polya's iteration and with 3 uncertain parameters in $d_1=d_2=4$ Polya's iterations. \vspace*{-0.17in}

\section{Conclusion}

In this paper, we have presented a cluster-computing and supercomputing approach to stability analysis of large-scale linear systems of the form $\dot x(t) = A(\alpha)x(t)$ where $A$ is polynomial, $\alpha \in \Delta_l \subset \R^l$ and $x \in \R^n$ and where $n\cong 100$ or $\alpha \cong 10$. The approach is based on mapping the structure of the LMI conditions associated with Polya's theorem to a decentralized computing environment. We have shown that for a sufficient number of processors, the proposed algorithm can solve the NP-hard robust stability problem with the same per-core computation cost as solving the Lyapunov inequality for a system with no parametric uncertainty. Theoretical and experimental results verify near-perfect scalability and speed-up for up to 200 processors.
Moreover, numerical examples demonstrate the ability of the algorithm to perform robust analysis of systems with 100+ states and several uncertain parameters using a simple 9-node Linux cluster computer.
We have also argued that our algorithms can also be extended to solve nonlinear stability analysis and robust controller synthesis problems, although this is left for future work.\vspace*{-0.1in}

\bibliographystyle{ieeetr} 
\bibliography{IEEE_Reza_paper_third_submission} \vspace*{-0.3in}
\begin{IEEEbiography}[{\includegraphics[width=1in,height=1.25in,clip,keepaspectratio]{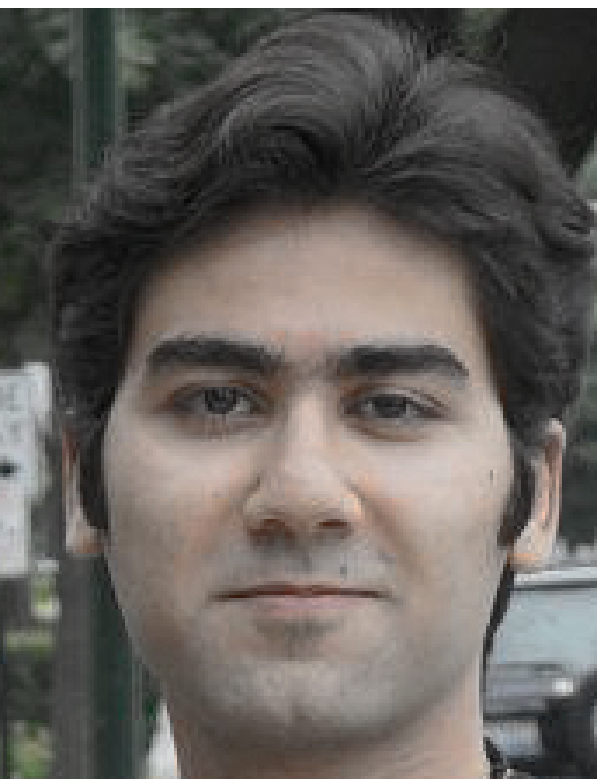}}]{Reza Kamyar}
received the B.S. and M.S in aerospace
engineering from Sharif University of Technology, Tehran, Iran in 2008, and 2010. He is currently a Ph.D student in the department of mechanical engineering of Arizona State University, Tempe, Arizona. He is a research assistant with Cybernetic Systems and Controls Laboratory (CSCL) in the School for Engineering of Matter, Transport and Energy (SEMTE) at Arizona State University. His research focuses on the development of decentralized algorithms applied to the problems of stability and control of large-scale complex systems. \vspace*{-0.3in}
\end{IEEEbiography}

\begin{IEEEbiography}
[{\includegraphics[width=1in,height=1.25in,clip,keepaspectratio]{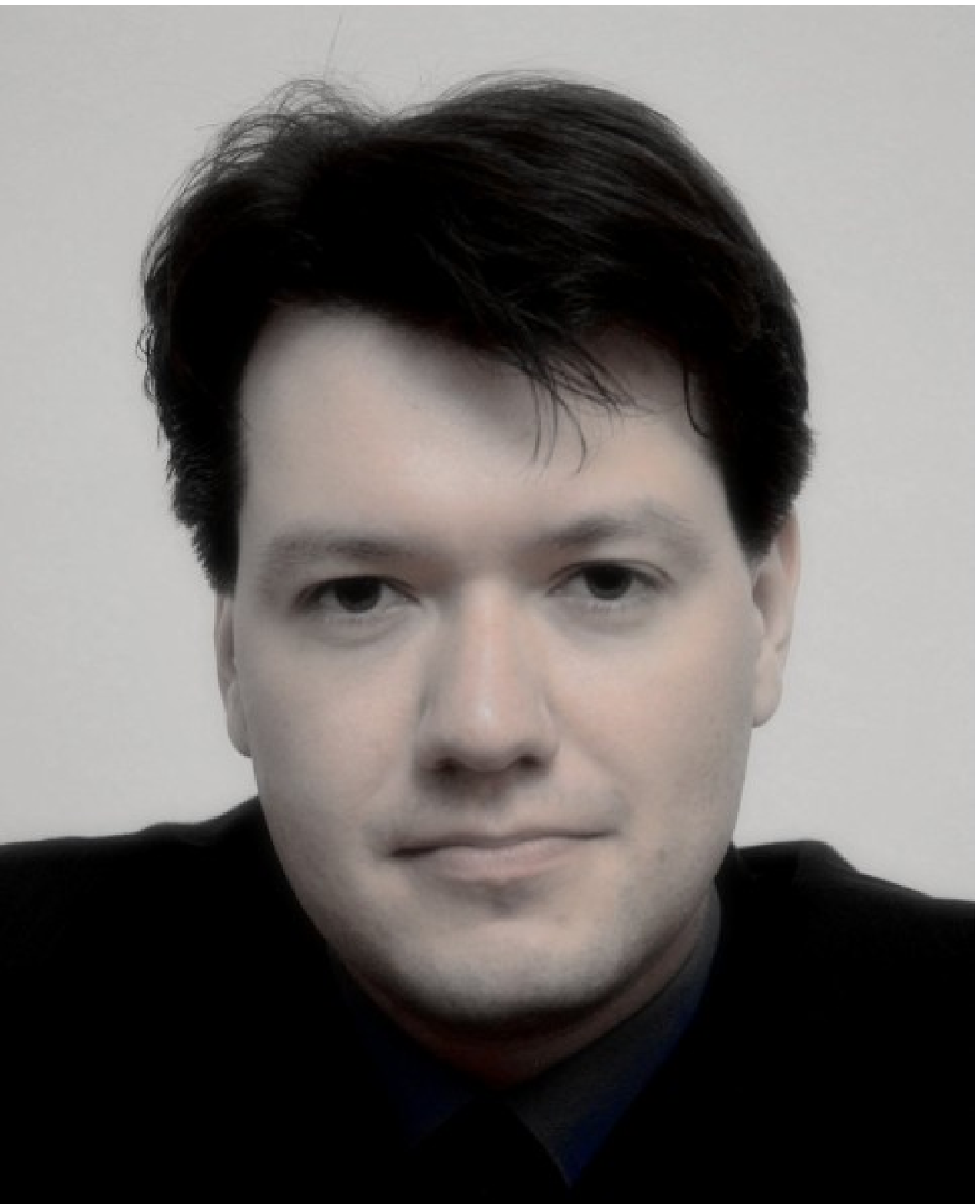}}]{Matthew M. Peet}
received the B.S. degrees in physics and in aerospace engineering from the University of Texas at Austin in 1999 and the M.S. and Ph.D. degrees in aeronautics and astronautics from Stanford University, Stanford, CA, in 2001 and 2006, respectively.
He was a Postdoctoral Fellow at the National Institute for Research in Computer Science and Control (INRIA), Paris, France, from 2006 to 2008, where he worked in the SISYPHE and BANG groups. He was an Assistant Professor of Aerospace Engineering in the Mechanical, Materials, and Aerospace Engineering Department, Illinois Institute of Technology, Chicago, from 2008 to 2012. He is currently an Assistant Professor of Aerospace Engineering in the School for the Engineering of Matter, Transport, and Energy at Arizona State University, Tempe, and Director of the Cybernetic Systems and Controls Laboratory. His research interests are in the role of computation as it is applied to the understanding and control of complex and large-scale systems. Applications include fusion energy and immunology. Dr. Peet received an NSF CAREER award in 2011. \vspace*{-0.3in}
\end{IEEEbiography}

\begin{IEEEbiography}
[{\includegraphics[width=1in,height=1.25in,clip,keepaspectratio]{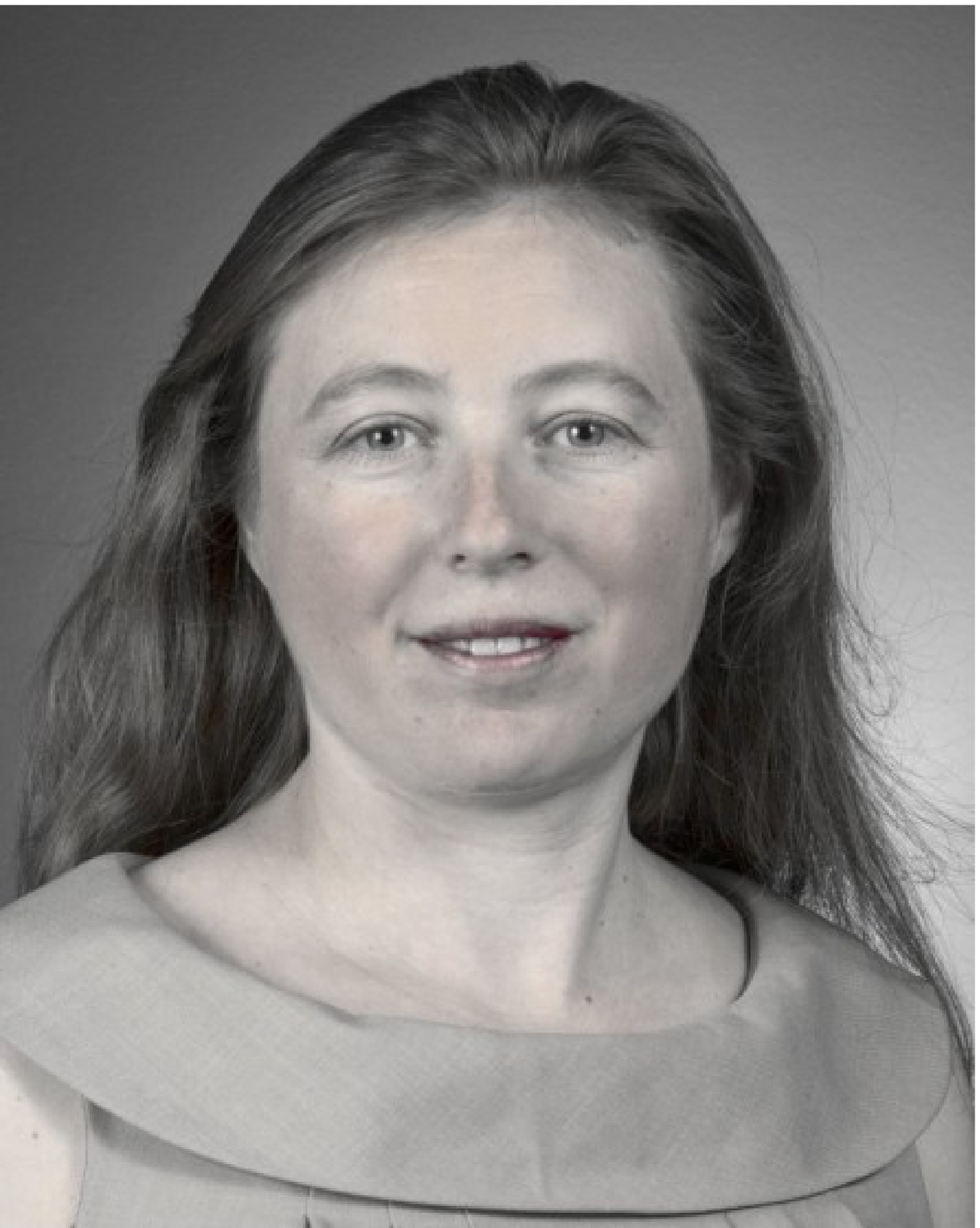}}]{Yulia Peet}
is an Assistant Professor of Mechanical and Aerospace Engineering at the School for Engineering of Matter, Transport and Energy at Arizona State University.  Her Ph.D. degree is in Aeronautics and Astronautics from Stanford (2006), M.S. in Aerospace Engineering (1999) and B.S. in Applied Mathematics and Physics (1997) from Moscow Institute of Physics and Technology in Russia. Her previous appointments include a postdoctoral position at the University of Pierre and Marie Curie in Paris in 2006-2008, and a dual appointment as an NSF research and teaching fellow at Northwestern University and assistant computational scientist at the Mathematics and Computer Science Division at Argonne National Laboratory in 2009-2012. 
\end{IEEEbiography}

\end{document}